\documentclass{amsart}
\usepackage{graphicx} 
\usepackage[usenames,dvipsnames]{color}
\usepackage{amsthm,amsfonts,amssymb,amsmath,amsxtra}
\usepackage[all]{xy}
\SelectTips{cm}{}
\usepackage{xr-hyper}

\usepackage[
backend=biber,
style=alphabetic,
doi = false
]{biblatex}

\usepackage[colorlinks=
 citecolor=Red,
  linkcolor=Black,
   urlcolor=Blue]{hyperref}
\usepackage{cleveref}

\usepackage[margin=1.25in]{geometry}
\usepackage{mathrsfs}

\usepackage{indentfirst}
\usepackage[symbol]{footmisc}

\RequirePackage{xspace}
\RequirePackage{etoolbox}
\RequirePackage{varwidth}
\RequirePackage{enumitem}
\RequirePackage{tensor}
\RequirePackage{mathtools}
\RequirePackage{longtable}
\RequirePackage{multirow}

\setitemize[0]{leftmargin=*,itemsep=\the\smallskipamount}
\setenumerate[0]{leftmargin=*,itemsep=\the\smallskipamount}

\newtheorem{theorem}{Theorem}
\newtheorem{proposition}[theorem]{Proposition}
\newtheorem{lemma}[theorem]{Lemma}

\newtheorem{corollary}[theorem]{Corollary}

\theoremstyle{definition}
\newtheorem{definition}[theorem]{Definition}

\newtheorem{remark}[theorem]{Remark}
\newtheorem{question}[theorem]{Question}
\newtheorem*{theorem*}{Theorem}

\numberwithin{equation}{section}
\numberwithin{theorem}{section}

\title{Theta Positivity in Lagrangian Grassmannian}
\author{Kaitao Xie}
\address{Department of Mathematics, The University of Hong Kong, Pokfulam, Hong Kong}
\email{kaitaoxie@connect.hku.hk}
\date{}

\begin{document}

\maketitle

\begin{abstract}
We study the theta nonnegative part of Lagrangian Grassmannian. We show that it admits an orbital decomposition and is homeomorphic to a closed ball. We compare it with other positive structures. We show that it contains several totally nonnegative parts of Lagrangian Grassmannian subject to certain choices of pinnings and agrees with the nonnegativity of the generalized Pl\"ucker coordinates.
\end{abstract}

\section{Introduction}
The theory of total positivity in split real reductive groups was introduced by Lusztig \cite{L94}. It has found connections to diverse research areas, including, among others, cluster algebras by Fomin and Zelevinsky \cite{FZ02}, higher Teichm\"uller theory by Fock and Goncharov \cite{FG06}, and theoretical physics by Arkani-Hamed and Trnka \cite{AT14}. Motivated by the application of total positivity in higher Teichm\"uller theory, Guichard and Wienhard \cite{GW18} introduced the notion of theta positivity. Theta positivity generalizes total positivity to a certain class of real reductive groups, not necessarily split, and provides a systematic way to understand the role of positivity in higher Teichm\"uller theory.

Besides the notion of total positivity in the reductive groups, Lusztig \cite{L94} also introduced the totally nonnegative part of the partial flag varieties with respect to a choice of pinning. It admits a representation-theoretical interpretation via the theory of canonical basis. He gave a decomposition of the totally nonnegative part of a partial flag variety by intersecting the projected Richardson varieties \cite[Remark 8.15]{L94}. He conjectured and Rietsch \cite{R97} showed that this decomposition is a cellular decomposition. 

Understanding the topology of the totally nonnegative part of a partial flag variety along with the cellular decomposition has been a central theme of the study of total positivity. Recent years have witnessed great development in this area. Marsh and Rietsch \cite[Theorem 11.3]{MR04} gave a parameterization of each cell. Rietsch \cite{R06} determined the closure relations of those cells. Williams \cite{Wil07} showed that the face poset is graded, shellable, and thin, and hence by a result of Bj\"orner \cite{Bjo84} is the face poset of a regular CW complex. She conjectured that the totally nonnegative part of a partial flag variety is a regular CW complex. Postnikov, Speyer, and Williams \cite{PSW09} showed that the totally nonnegative part of Grassmannian is a CW complex. Their result was generalized by Rietsch and Williams \cite{RW08} to general partial flag varieties. Hersh \cite{Her14} obtained the regularity result for the link of identity in the totally nonnegative unipotent monoid. Galashin, Karp, and Lam \cite{GKL19,GKL22a} showed that the totally nonnegative part of a partial flag variety is homeomorphic to a closed ball. Williams' conjecture was finally proven by Galashin, Karp, and Lam \cite{GKL22b} for finite types, and by Bao and He \cite{BH24} for Kac-Moody types. 

Another interesting positive structure on partial flag varieties is the positive (nonnegative) points of the generalized Pl\"ucker coordinates. It leads to connections with various topics, such as tropical geometry by Speyer and Williams \cite{SW05}, loop amplituhedra by Arkani-Hamed and Trnka \cite{AT14}, also by Bai, He and Lam \cite{BHL16}, and so on. In the case of Grassmannian, this is Postnikov's \cite{P06} definition of totally positive (nonnegative) Grassmannian. He used the matroid strata of Grassmannian to give a cellular decomposition of the totally nonnegative Grassmannian. Rietsch showed in an unpublished note that, Lusztig's and Postnikov's definitions and cellular decompositions of totally nonnegative Grassmannian agree. Subsequent proofs were given by Talaska and Williams \cite{TW13}, Lam \cite{Lam16}, and Lusztig \cite{L19}. The comparison was carried on by Bloch and Karp \cite{BK23b} for partial flag varieties of type $A$. They showed that the two notions of positivity agree if and only if the partial flag variety is consecutive. They also showed that the generalized Pl\"ucker nonnegative part of a partial flag variety of type $A$ is homeomorphic to a closed ball, no matter whether it agrees with the totally nonnegative part or not \cite[Corollary 6.16]{BK23a}. 

The goal of this paper is to study the topology of the theta nonnegative part of Lagrangian Grassmannian, to find analogs of the results on the totally nonnegative part of partial flag varieties, as well as to compare the three positive structures. The main results are summarized as follows.

\begin{theorem*}
Let $\mathcal P$ denote the Lagrangian Grassmannian. Let $\mathcal P_{>0}$ and $\mathcal P_{\geq 0}$ be the theta positive part and theta nonnegative part of $\mathcal P$, respectively.
\begin{enumerate}
    \item The space $\mathcal P_{\geq 0}$ decomposes into finitely many $GL_n(\mathbb R)^+$-orbits, and each one is a connected component of a sub-variety of $\mathcal P(\mathbb R)$, which can be regarded as a variation of the open Richardson varieties.
    \item The space $\mathcal P_{\geq 0}$ is homeomorphic to a closed ball, and $\mathcal P_{\geq 0}\setminus\mathcal P_{>0}$ is homeomorphic to a sphere.
    \item There is a family of pinnings, such that the totally positive parts of Lagrangian Grassmannian subject to these pinnings lie in $\mathcal P_{>0}$. In particular, $\mathcal P_{\geq0}$ is not the union of several totally nonnegative parts of $\mathcal P$ subject to a certain family of pinnings, except the Lagrangian Grassmannian of a $4$-dimensional symplectic vector space.
    \item The notions of theta positivity and generalized Pl\"ucker positivity agree in $\mathcal P$.
\end{enumerate}
\end{theorem*}

The last two statements together answer a question in \cite{BK23b} about comparing the total positivity and generalized Pl\"ucker positivity in the case of Lagrangian Grassmannian. 

However, none of the open Richardson varieties, the projected Richardson varieties, or the general version of matroid strata (called the Gelfand-Serganova strata) give rise to a cellular decomposition of $\mathcal P_{\geq 0}$. It would be interesting to find a nice cellular decomposition for $\mathcal P_{\geq 0}$. 

The theta positive part of Lagrangian Grassmannian also appears in the study of Gaussoid \cite{BDKS19} from an algebraic statistics point of view. While the results there do not play a role in this paper, it would be interesting to explore further connections.

\subsection*{Outline} In \Cref{background}, we review the background on total positivity, generalized Pl\"ucker coordinates, and theta positivity. In \Cref{main}, we define the object of major interest in this paper and present our main results. In \Cref{combinatorics}, we collect the relevant combinatorial results. In \Cref{decomposition}, we use the results from \Cref{combinatorics} to give an orbital decomposition of the theta nonnegative part of Lagrangian Grassmannian, and we also obtain a cellular decomposition of the theta nonnegative unipotent monoid of the symplectic group. In \Cref{topology}, we show that the theta nonnegative part of Lagrangian Grassmannian is a closed ball by employing the techniques of contractive flows. In \Cref{comparison}, we compare the object studied here with other positive structures. 

\subsection*{Acknowledgement}
I thank Xuhua He for many inspiring discussions. I thank Thomas Lam and Weinan Zhang for valuable comments.

\section{Background} \label{background}

\subsection{Total positivity} \label{total_positivity}
Let $G$ be a reductive, connected, linear algebraic group over $\mathbb C$ split over $\mathbb R$. We fix a maximal torus $T$, and a pair of opposite Borel subgroups $B^\pm$ with unipotent radicals $U^\pm$, such that $T = B^+\cap B^-$. Denote by $X(T)$ and $Y(T)$ the sets of weights and co-weights of $T$, respectively. Let $W:= N_G(T)/T$ be the Weyl group of $G$. For each $w\in W$, denote by $\dot w$ one of its lift in $G$. 

Let $I$ be the indexing set of simple roots. For each $i\in I$, denote by $\alpha_i:T\to\mathbb R$ the simple root, by $\alpha_i^\vee$ the simple coroot, and by $U_i^\pm$ the positive/negative simple root subgroups of $U^\pm$. We fix a set of isomorphisms
$$x_i: \mathbb C\to U_i^+, y_i: \mathbb C\to U_i^-, i\in I,$$
such that the assignment
$$\begin{bmatrix} 1 & a \\ 0 & 1 \end{bmatrix}\to x_i(a), \begin{bmatrix} 1 & 0 \\ b & 1 \end{bmatrix}\to y_i(a),
\begin{bmatrix} c & 0 \\ 0 & c^{-1} \end{bmatrix}\to \alpha^\vee_i(c), a,b\in\mathbb C, c\in\mathbb C^*,$$
defines a homomorphism $SL_2\to G$.

Let $\mathcal B := G/B^+$ be the full flag variety of $G$. For $v,w\in W$ such that $v\leq w$ in the Bruhat order, we have the open Richardson variety $\mathcal B_{v,w}:=B^+\dot w\cdot B^+\cap B^-\dot v\cdot B^+\subseteq\mathcal B$. For $J\subseteq I$, we have a standard parabolic subgroup $P_J^+$ containing $B^+$, and $W_J$, the corresponding parabolic subgroup of $W$. Denote by $W^J$ the set of minimal representatives of $W/W_J$. We have the partial flag variety $\mathcal P:=G/P_J^+$, the natural projection $\pi:G/B^+\to G/P_J^+$, and the projected Richardson varieties $\mathcal P_{v,w}:=\pi(\mathcal B_{v,w})$ for $v\leq w$, $v\in W$ and $w\in W^J$. The projected Richardson varieties form a stratification of $\mathcal P$. General results about the projected Richardson varieties can be found in \cite{KLS14}.

The data $\mathbf P = (T,B^+,B^-,x_i,y_i; i\in I)$ is called a pinning. Let $T_{>0}$ be the subgroup of $T$ generated by $\chi(c)$ for $\chi\in Y(T)$ and $c\in\mathbb R_{>0}$. Take two reduced words $\mathbf i = (i_1,i_2,\cdots,i_k)$ and $\mathbf i^\prime = (i^\prime_1,i^\prime_2,\cdots,i^\prime_k)$ for the longest element in $W$. Define
$$U_{(\mathbf P,>0)}^+ = \{x_{i_1}(a_1)x_{i_2}(a_2)\cdots x_{i_k}(a_k)|a_j\in\mathbb R_{>0}\}\simeq\mathbb R_{>0}^k,$$
$$U_{(\mathbf P,>0)}^- = \{y_{i^\prime_1}(b_1)y_{i^\prime_2}(b_2)\cdots y_{i^\prime_k}(b_k)|b_j\in\mathbb R_{>0}\}\simeq\mathbb R_{>0}^k.$$

The semi-groups $U^\pm_{(\mathbf P,>0)}$ are independent of $\mathbf i$ and $\mathbf i^\prime$ chosen \cite[Lemma 2.10]{L94}.

\begin{definition}
The totally nonnegative part $G_{(\mathbf P,\geq 0)}$ of $G$ is the submonoid of $G$ generated by all $x_i(a)$, $y_i(b)$, $\chi(c)$ for $a,b,c\in\mathbb R_{>0}$, $i\in I$ and $\chi\in Y(T)$. The totally positive part of $G$ is $G_{(\mathbf P,>0)} := U^-_{(\mathbf P,>0)}T_{>0}U^+_{(\mathbf P,>0)}$. 

\end{definition}

The triple factorization of $G_{(\mathbf P,>0)}$ can be taken in any order \cite[\S1.3,\S2.12]{L94}.

\begin{definition}
The totally positive part $\mathcal P_{(\mathbf P,>0)}$ and the totally nonnegative part $\mathcal P_{(\mathbf P,\geq 0)}$ of $\mathcal P$ are defined by

$$\mathcal P_{(\mathbf P,>0)} := G_{(\mathbf P,>0)}\cdot P_J^+\subseteq\mathcal P(\mathbb R),\ \mathcal P_{(\mathbf P,\geq 0)} := \overline{G_{(\mathbf P,\geq 0)}\cdot P_J^+}\subseteq\mathcal P(\mathbb R).$$ 

\end{definition}

Let $\mathcal P_{v,w,(\mathbf P,>0)} := \mathcal P_{(\mathbf P,\geq 0)}\cap\mathcal P_{v,w}(\mathbb R)$. Then $\mathcal P_{(\mathbf P,\geq0)} = \bigsqcup\mathcal P_{v,w,(\mathbf P,>0)}$. As mentioned in the introduction, $\mathcal P_{(\mathbf P,\geq0)}$ admits a regular CW complex structure. More explicitly, we have

\begin{enumerate}
    \item $\mathcal P_{v,w,(\mathbf P,>0)}$ is a semi-algebraic cell and a connected component of $\mathcal P_{v,w}(\mathbb R)$ \cite{R97,MR04};
    \item the Hausdorff closure of $\mathcal P_{v,w,(\mathbf P,>0)}$ is a union of $\mathcal P_{v^\prime,w^\prime,(\mathbf P,>0)}$ for some $(v^\prime,w^\prime)$ \cite{R06};
    \item the decomposition $\overline{\mathcal P_{v,w,(\mathbf P,>0)}} = \bigsqcup \mathcal P_{v^\prime,w^\prime,(\mathbf P,>0)}$ is a regular CW complex \cite{GKL22b,BH24}.
\end{enumerate}

In particular, $\mathcal P_{(\mathbf P,\geq 0)}$ is homeomorphic to a closed ball.

\subsection{Generalized Pl\"ucker coordinates} \label{generalized_plucker}
The notions of generalized Pl\"ucker coordinates and strata for partial flag varieties of general Lie types were defined in \cite{GS87}. The geometry of the Gelfand-Serganova strata is related to the combinatorial theory of Coxeter matroids developed in \cite{BGW03}.

Keep the notations from \Cref{total_positivity}. Denote by $\omega_i\in X(T)$ the fundamental weight corresponding to $\alpha^\vee_i$, for $i\in I$. For $\lambda\in X^+(T)$ a dominant weight, let $V_\lambda$ be the irreducible representation with highest weight $\lambda$ and highest weight vector $\eta_\lambda$. The weights in the Weyl group orbit of $\lambda$ are called the extremal weights.

For $J\subseteq I$, let $\rho_J := \sum_{i\in I\setminus J}\omega_i$. Then we have the embedding $G/P_J^+\to\mathbb P(V_{\rho_J})$ given by $gP_J^+\to g\cdot \eta_{\rho_J}$. Let $\mathrm A$ be the set of weights appearing in $V_{\rho_J}$. We pick a set of weight vectors $\{\eta_\mu|\mu\in\mathrm A\}$. Every point $x\in G/P_J^+$ uniquely determines, up to a common constant $d$, a collection of numbers $p^\mu(x)$, $\mu\in\mathrm A$, where 
$$x = d\sum_{\mu\in\mathrm A}p^\mu(x)\eta_\mu.$$ 

The collection of numbers associated with extremal weights $\{p^\mu(x)|\mu\in W\cdot\rho_J\}$, defined up to a common scalar, are called generalized Pl\"ucker coordiates of $x$. The list of $x$ is defined to be
$$L_x := \{\mu\in W\cdot\rho_J|p^\mu(x)\neq 0\}.$$
A Gelfand-Serganova stratum of $\mathcal P$ is defined to be the subset of points with the same list. 

For $i\in I$, set $\dot s_i = x_i(-1)y_i(1)x_i(-1)$, and $\dot w = \dot w_1\dot w_2$, if $w = w_1w_2$ and $l(w) = l(w_1)+l(w_2)$, where $l:W\to\mathbb N$ is the length function. In other words, we pick a lift $\dot w$ for each $w\in W$ in a certain way. We obtain a well-defined set of extremal weight vectors $\{\eta_{w\lambda} := \dot w\eta_\lambda|w\in W\}$. From now on, we assume that the extremal weights are chosen in this way. 

We define the generalized Pl\"ucker positive (nonnegative) part of $\mathcal P$ to be the subset $\mathcal P_{\triangle>0}$ ($\mathcal P_{\triangle\geq0}$, respectively) consisting of points $x\in\mathcal P$ whose generalized Pl\"ucker coordinates are all positive (nonnegative, respectively) up to a common scalar. Note that the definition here implicitly involves a pinning $\mathbf P$, as we need it to choose the set $\{\dot w|w\in W\}$ systematically. 

We review some previous studies on this subject.

\begin{enumerate}
    \item When $G = SL_n$ and $P_J$ is a maximal parabolic subgroup, the partial flag variety $SL_n/P_J$ is a Grassmannian. The generalized Pl\"ucker coordinates are classical Pl\"ucker coordinates, and the Gelfand-Serganova strata are the matroid strata. The generalized Pl\"ucker nonnegative part $(SL_n/P_J)_{\triangle\geq0}$ is Postnikov's definition of totally nonnegative Grassmannian, and the matroid strata give a cellular decomposition \cite{P06}. The definition and cellular decomposition coincide with Lusztig's total positivity (see \cite{TW13,Lam16,L19}).

    \item  When $G = SL_n$ and $P_J$ is a general parabolic subgroup, $(SL_n/P_J)_{\triangle\geq0}$ was studied in \cite{BK23a,BK23b}, with a different but equivalent definition. See \cite[Remark 5.13]{BK23b}. It was shown to be a closed ball for arbitrary $J$ \cite[Corollary 6.16]{BK23a} and coincides with total positivity if and only if $I\setminus J$ is a connected Dynkin diagram \cite{BK23b}.

    \item For general Lie types, Tsukerman and Williams \cite[Theorem 7.1]{TW15} showed that one totally positive cell of $\mathcal P_{(\mathbf P,\geq 0)}$ lies in one Gelfand-Serganova stratum.
\end{enumerate}

It is not immediate from the definition that the Hausdorff closure of $\mathcal P_{\triangle>0}$ is $\mathcal P_{\triangle\geq0}$. The rest of this subsection is devoted to showing this, starting by recalling the definition of generalized minors \cite[Definition 6.2]{MR04}.

\begin{definition}
For $\lambda\in X^+(T)$, $v,w\in W$, $g\in G$, define $\triangle^{v\lambda}_{w\lambda}(g)$ to be the coefficients of $\eta_{v\lambda}$ in $g\eta_{w\lambda}$ inside $V_\lambda$. 
\end{definition}

\begin{lemma}
For $g\in G_{(\mathbf P,>0)}$, we have $\triangle^{v\lambda}_{w\lambda}(g)>0$ for every $v,w\in W$ and $\lambda\in X^+(T)$.

\begin{proof}
If $\lambda = \mu_1+\mu_2$ with $\mu_1,\mu_2\in X^+(T)$, then there is a $G$-equivariant embedding $V_\lambda\to V_{\mu_1}\otimes V_{\mu_2}$ by $\eta_\lambda\to \eta_{\mu_1}\otimes\eta_{\mu_2}$. Since every extremal weight space is of dimension $1$, we can conclude that 
$$\triangle^{v\lambda}_{w\lambda}(g) = \triangle^{v\mu_1}_{w\mu_1}(g)\triangle^{v\mu_2}_{w\mu_2}(g).$$

Now we only need to consider the case where $\lambda = \omega_i$ for some $i\in I$. This was proven by Fomin and Zelevinsky \cite[Theorem 1.11-1.12]{FZ99}. 
\end{proof}

\end{lemma}

\begin{proposition} \label{plucker_closure}
The Hausdorff closure of $\mathcal P_{\triangle>0}$ is $\mathcal P_{\triangle\geq0}$.

\begin{proof}
By definition, the space $\mathcal P_{\triangle\geq0}$ is closed and contains $\mathcal P_{\triangle>0}$, so we have $\overline{\mathcal P_{\triangle>0}}\subseteq\mathcal P_{\triangle\geq0}$. By \cite[Proposition 5.9(c)]{L94}, there is a continuous function $f:\mathbb R_{\geq 0}\to G$ such that $f(0) = e$ and $f(t)\in G_{(\mathbf P,>0)}$ for $t>0$. For $x\in\mathcal P_{\triangle\geq0}$, by the lemma above we have $f(t)x\in\mathcal P_{\triangle>0}$ for $t>0$. Since $f(0)x = x$, the conclusion follows.
\end{proof}
\end{proposition}

\begin{remark}
We see that $\mathcal P_{(\mathbf P,>0)}\subseteq\mathcal P_{\triangle>0}$ and $\mathcal P_{(\mathbf P,\geq0)}\subseteq\mathcal P_{\triangle\geq0}$. In particular, they are not empty.
\end{remark}

As we will show in \Cref{comparison}, the theta positivity of Lagrangian Grassmannian studied in this paper agrees with the generalized Pl\"ucker positivity.

\subsection{Theta positivity\protect\footnote{The name ``theta" refers to a subset $\Theta$ of simple roots \cite{GW24}. To keep our notations coherent, we prefer to use $J$ to indicate the complement of $\Theta$.}}

We again keep the notations in \Cref{total_positivity}, except in this section we take $G$ a real reductive Lie group without the assumption that $G$ is split over $\mathbb R$. For $J\subseteq I$, We have a pair of parabolic subgroups $P^+_J$ and $P_J^-$. Let $L_J: = P_J^+\cap P_J^-$ be the Levi subgroup, with identity component $L_J^\circ$, and $U_J^\pm$ be the unipotent radical of $P_J^\pm$. Let $\mathfrak{u}_J^\pm$ be the Lie algebra of $U_J^\pm$. Then $\mathfrak u_J^\pm$ is a representation of $L_J$. This representation decomposes into irreducible component, and for each $i\in I\setminus J$, we have an irreducible component $\mathfrak u_{J,i}^\pm$, that is,
$\mathfrak u_J^\pm = \bigoplus_{i\in I\setminus J}\mathfrak u_{J,i}^\pm$. We say that $G$ admits a theta positive structure if 
\begin{itemize}
    \item There is a subset $J\subseteq I$, such that for each $i\in I\setminus J$, there is a closed, acute, convex cone $c^\pm_i\subseteq \mathfrak u_{J,i}^\pm$ that is invariant under $L^\circ_J$-action. Here, acute means that it contains no line passing through the origin.
\end{itemize}

Guichard and Wienhard \cite{GW24} showed that there are four families of groups admitting theta positive structures:
\begin{enumerate}
    \item $G$ is split over $\mathbb R$ and $J = \emptyset$;
    \item $G$ is Hermitian of tube types, and $J$ consists of the indices of all short simple roots;
    \item $G$ is locally isomorphic to $SO(p+1,p+k)$, $p>1$, $k>1$, and $J$ consists of the indices of all short simple roots;
    \item $G$ is the real form of $F_4$, $E_6$, $E_7$ or $E_8$ whose reduced root system is of type $F_4$, and $J$ consists of the indices of all short simple roots.
\end{enumerate}

They then introduced

\begin{definition}
The theta nonnegative monoid $G_{J,\geq 0}$ is the submonoid of $G$ generated by $L_J^\circ$ and $\exp(x)$ for $x\in c_i^\pm$, $i\in I\setminus J$. 
\end{definition}
 
For the first family, where $G$ is split over $\mathbb R$, the theta nonnegative monoid recovers Lusztig's totally nonnegative monoid. For the second family, the theta nonnegative monoid recovers the contraction semi-group, studied in, for example, \cite{Kou95}.

\begin{remark}
Each one of the split real Lie groups of types $C_n$, $B_n$, and $F_4$ has two theta positive structures. The one with $J = \emptyset$ is total positivity.
\end{remark}

Following the study of total positivity in partial flag varieties, a natural question to ask is

\begin{question}
What can we say about the topology of $\overline{G_{J,\geq 0}\cdot P_J^+}\subseteq\mathcal P$?
\end{question}

In this paper, we study this question in the case where $G = Sp_{2n}(\mathbb R)$ and $J$ consists of the indices of all short simple roots.

\section{Main Results} \label{main}
Consider a linear space $\mathbb R^{2n}$ with a symplectic form $\Omega$ given by the matrix $\begin{bmatrix}
    0 & I_n \\
    -I_n & 0
\end{bmatrix}$, where $I_n$ is the $n$ by $n$ identity matrix. The symplectic group $G = Sp_{2n}(\mathbb R)$ is identified as the matrices group 
$$G = \left\{\begin{bmatrix}A&B\\ C&D\end{bmatrix}\in GL_{2n}(\mathbb R)\bigg| A^tD-B^tC = I_n, A^tC = C^tA, B^tD = D^tB\right\}.$$
Fix the maximal torus $T = \{\text{diag}(a_1,a_2,\cdots,a_n,a_1^{-1},a_2^{-1},\cdots,a_n^{-1})|a_i\in\mathbb R^*,1\leq i\leq n\}$
and a pair of opposite Borel subgroups
$$B^+=\left\{\begin{bmatrix}A&B\\ 0&(A^t)^{-1}\end{bmatrix}\in Sp_{2n}(\mathbb R)\bigg|A\text{ is upper triangular} \right\},$$ 

$$B^-=\left\{\begin{bmatrix}A&0\\ C&(A^t)^{-1}\end{bmatrix}\in Sp_{2n}(\mathbb R)\bigg|A\text{ is upper triangular} \right\}.$$

Denote by $E_{i,j}$ the $n$ by $n$ matrix whose entries are all $0$ except that the $(i,j)$-entry is $1$. We choose a set of Chevalley generators $e_i = E_{i,i+1}-E_{n+i+1,i}$ for $1\leq i\leq n-1$, $e_n = E_{n,2n}$ and $f_i = e_i^t$. In this way, we obtain a pinning $(T,B^+,B^-,x_i,y_i;1\leq i\leq n)$, where $x_i:\mathbb R\to G$ and $y_i:\mathbb R\to G$ are given by $x_i(a) = \exp(ae_i)$ and $y_i(a) = \exp(af_i)$. We pick a lift $\dot w$ for each $w\in W$ by letting $\dot{s}_i = x_i(-1)y_i(1)x_i(-1)$ and $\dot w = \dot w_1\dot w_2$, if $w = w_1w_2$ and $l(w) = l(w_1)+l(w_2)$. 

Let $J = \{1,2,\cdots, n-1\}$. Consider a pair of parabolic subgroups and their Levi subgroup $$P_J^+ = \left\{\begin{bmatrix}A&B\\ 0&A^{t,-1}\end{bmatrix}\bigg|AB^t=BA^t\right\},\  P_J^- = \left\{\begin{bmatrix}A&0\\ C&A^{t,-1}\end{bmatrix}\bigg| A^tC = C^tA\right\},\ L_J = P_J^+\cap P_J^-\simeq GL_n(\mathbb R).$$ The variety of major interest is the partial flag variety $\mathcal P:=G/P_J^+$. By \cite[6.2.11.(3)]{Spr98}, $\mathcal P$ is isomorphic to the Lagrangian Grassmannian of the real symplectic space $(\mathbb R^{2n},\Omega)$. As a result, 
$$\mathcal P = \left\{\begin{bmatrix}A\\ C\end{bmatrix}\in Gr(n,2n)(\mathbb R)\bigg| A,C\in \text{Mat}_{n\times n}(\mathbb R), A^tC = C^tA\right\}.$$
Note that $(A,C)^t = (A^\prime,C^\prime)^t$ in $\mathcal P$ if and only if there is an $n\times n$ invertible matrix $g$ such that $A^\prime = Ag$ and $C^\prime = Cg$. 

Following \cite{GW24}, we define the theta positive unipotent semi-groups
$$U^-_{J,>0} := \left\{\begin{bmatrix} I_n&0\\C&I_n
\end{bmatrix}\bigg|C\text{ is symmetric and positive definite}\right\},$$
$$U^+_{J,>0} := \left\{\begin{bmatrix} I_n&B\\0&I_n
\end{bmatrix}\bigg|B\text{ is symmetric and positive definite}\right\},$$
and the theta nonnegative unipotent monoids to be the respective closure 
$$U^-_{J,\geq 0} := \overline{U^-_{J,>0}}=\left\{\begin{bmatrix} I_n&0\\C&I_n
\end{bmatrix}\bigg|C\text{ is symmetric and positive semi-definite}\right\},$$ 
$$U^+_{J,\geq 0} := \overline{U^+_{J,>0}}=\left\{\begin{bmatrix} I_n&B\\0&I_n
\end{bmatrix}\bigg|B\text{ is symmetric and positive semi-definite}\right\}.$$
The theta nonnegative monoid $G_{J,\geq 0}$ is defined to be the submonoid of $G$ generated by $U^-_{J,\geq 0}$, $U^+_{J,\geq 0}$, and $L_J^\circ$, the identity component of $L_J$. 

In \cite{Kou95}, it was shown that $G_{J,\geq 0}$ has a triple factorization $G_{J,\geq 0} = U_{J,\geq 0}^-L_J^\circ U_{J,\geq 0}^+$, and
$$G_{J,\geq 0} = \left\{\begin{bmatrix}A&B\\ C&D\end{bmatrix}\in Sp_{2n}(\mathbb R)\bigg| D\text{ is invertible}, CD^t, D^tB\text{ are positive semi-definite}\right\}.$$

We define the theta positive part of $\mathcal P$ by
$$\mathcal P_{>0} := U_{J,>0}^-\cdot P^+_J = \left\{\begin{bmatrix}I\\ C\end{bmatrix}\in\mathcal P\bigg| C\text{ is positive definite}\right\}.$$

It is clear that $\mathcal P_{>0}$ is open in $\mathcal P$. In \cite[Proposition 13.1]{GW24}, it was shown that $$\mathcal P_{>0} = U_{J,>0}^-\cdot P_J^+ = U_{J,>0}^+\cdot P_J^-,$$ and that $\mathcal P_{>0}$ is a connected component of $P^-\cdot P_J^+\cap P_J^+\dot w_0\cdot P_J^+$, where $w_0$ is the longest element in $W$.

We first define the theta nonnegative part by
$$\mathcal P_{\geq 0} := \left\{\begin{bmatrix}A\\ C\end{bmatrix}\in\mathcal P\bigg| A^tC = C^tA\text{ is positive semi-definite}\right\}.$$

Note that this condition is independent of the representative chosen. If $A^\prime = Ag$ and $C^\prime = Cg$ for some invertible matrix $g$, then $A^{\prime t}C^\prime = g^tA^tCg$, and $A^{\prime t}C^\prime$ will be positive semi-definite if $A^tC$ is. 

By the triple decomposition, we have that $ \overline{U_{J,>0}^-\cdot P^+_J} = \overline{G_{J,\geq 0}\cdot P^+_J}$. However, it is a nontrivial fact that the Hausdorff closure of $\mathcal P_{>0}$ is $\mathcal P_{\geq 0}$. We will show this in \Cref{closure}.

\begin{remark}
One can show that a matrix $\begin{bmatrix} A&B \\ C&D\end{bmatrix}\in G_{J,\geq 0}$ guarantees that $A^tC$ is positive semi-definite. However, for such a matrix in $G_{J,\geq0}$, $A$ is invertible, so $\mathcal P_{\geq 0}$ is not a $G_{J,\geq 0}$-orbit.
\end{remark}

In this paper, we study the topology of $\mathcal P_{>0}$ and $\mathcal P_{\geq 0}$. We now formulate the results as the following theorem.

\begin{theorem}
 The spaces $\mathcal P_{>0}$ and $\mathcal P_{\geq 0}$ have the following properties.
\begin{enumerate}

    \item The orbits of $L^\circ_J$-action on $\mathcal P_{\geq 0}$ are $\mathcal P_{\geq 0}\cap P_J^-\dot v\cdot P_J^+\cap P_J^+\dot w\cdot P_J^+$ for various $v\leq w$, and each one is a connected component of the respective $P_J^-\dot v\cdot P_J^+\cap P_J^+\dot w\cdot P_J^+$.
    (\Cref{orbit} and \Cref{connected_component}). As an application, the Hausdorff closure of $\mathcal P_{>0}$ is $\mathcal P_{\geq 0}$ (\Cref{closure}).

    \item The space $\mathcal P_{\geq 0}$ is homeomorphic to a closed ball of dimension $\dfrac{n(n+1)}{2}$, and $\mathcal P_{\geq 0}\setminus\mathcal P_{>0}$ is homeomorphic to a sphere of dimension $\dfrac{n(n+1)}{2}-1$ (\Cref{ball}).

    \item There is a collection of pinnings $\mathrm C$ such that
    $$\bigsqcup_{\mathbf P\in\mathrm C}\mathcal P_{(\mathbf P,>0)}\subsetneqq\mathcal P_{>0},\  \bigsqcup_{\mathbf P\in\mathrm C}\mathcal P_{(\mathbf P,\geq 0)}\subseteq\mathcal P_{\geq0}.$$ The equality of nonnegative parts holds when $n=2$, but fails when $n\geq 3$ (\Cref{dense} and \Cref{TP}).

    \item With respect to the extremal weight vectors $\{\dot w\eta_{\omega_n}|w\in W^J\}$, we have that $\mathcal P_{>0} = \mathcal P_{\triangle>0}$ and $\mathcal P_{\geq0} = \mathcal P_{\triangle\geq0}$ (\Cref{Plucker}).

\end{enumerate}
\end{theorem}

We do not have a cellular decomposition for $\mathcal P_{\geq 0}$. Nevertheless, we have a cellular decomposition of $U_{J,\geq 0}^-$ by intersecting $B^+\dot wP_J^+$, which coincides with the Cholesky decomposition of symmetric positive semi-definite matrices (\Cref{cellular}).

\section{Weyl group of type $B/C$} \label{combinatorics}
 Let $W$ be the Weyl group of type $B_n/C_n$, with the set of simple reflections $I=\{s_1,s_2,\cdots,s_{n-1},t\}$ such that $s_is_{i+1}s_i = s_{i+1}s_is_{i+1}$ for $i = 1,\cdots,n-2$ and $s_{n-1}ts_{n-1}t = ts_{n-1}ts_{n-1}$. Let $J = \{s_1,s_2,\cdots,s_{n-1}\}$.

Consider the permutation group $\text{Perm}(\pm1,\pm2,\cdots,\pm n)$ on $\{\pm1,\pm2,\cdots,\pm n\}$. We use plain cycle $(i_1,\cdots,i_k)$ to denote the permutation $\sigma$ such that $\sigma(a) = a$ if $a\notin\{\pm i_1,\pm i_2,\cdots,\pm i_k\}$, $\sigma(\pm i_j) = \pm i_{j+1}$ for $j = 1,\cdots,k-1$ and $\sigma(\pm i_k) = \pm i_1$. We use the bar cycle $(\overline{i_1,\cdots,i_k})$ to denote the permutation $\sigma$ such that $\sigma(a) = a$ if $a\in\{\pm1,\pm2,\cdots,\pm n\}\setminus\{\pm i_1,\cdots,\pm i_k\}$, $\sigma(\pm i_j) = \pm i_{j+1}$ for $j = 1,\cdots,k-1$, $\sigma(i_k) = -i_1$ and $\sigma(-i_k) = i_1$. 

It is well-known that $W$ can be identified as $\{\sigma\in \text{Perm}(\pm1,\pm2,\cdots,\pm n)|\sigma(-i) = -\sigma(i), \leq i\leq n\}$ with $s_i = (i,i+1)$ and $t = (\overline{n})$.

\begin{remark}
The lift $\dot w$ considered in \Cref{main} is a signed permutation matrix on $\{1,2,\cdots,2n\}$. It shall be identified with the Weyl group element here via the bijection $f:\{1,2,\cdots,2n\}\to\{\pm 1,\pm 2,\cdots,\pm n\}$ defined by $f(k) = k$ if $k\leq n$ and $f(k) = n-k$ if $k>n$.
\end{remark}

\subsection{Minimal Representatives in Single and Double Cosets}
We identity $W/W_J$ with $2^{\{1,2,\cdots,n\}}$ by associating each $K\subseteq\{1,2,\cdots,n\}$ with a minimal representative in $W/W_J$.

For $K\subseteq\{1,2,\cdots,n\}$, define $w_K = \prod_{k\in K}^{\downarrow}s_ks_{k+1}\cdots s_{n-1}t$, where we adopt the convention that $s_n:= t$. Here the order of product is the decreasing order of $K$, and we use $\downarrow$ to emphasize it. In terms of cycles, $s_ks_{k+1}\cdots s_{n-1}t = (\overline{k,k+1,\cdots,n})$, and $(s_ks_{k+1}\cdots s_{n-1}t)^{-1} = (\overline{n,n-1,\cdots,k})$.

\begin{lemma} \label{left_multi}
Suppose that $i\in\{1,2,\cdots,n\}$. Then
\begin{enumerate}[label = (\alph*)]
    \item if $i=n\notin K$, then $tw_K = w_{K\cup\{n\}}$; if $i=n\in K$, then $tw_K = w_{K\setminus\{n\}}$;
    \item if $i\neq n$, $i\notin K$ and $i+1\in K$, then $s_iw_K = w_{K\cup\{i\}\setminus\{i+1\}}$;
    \item if $i\neq n$, $i\in K$ and $i+1\notin K$, then $s_iw_K = w_{K\cup\{i+1\}\setminus\{i\}}$;
    \item otherwise, $s_iw_K = w_Ks_j$ where $j = n-1-\#(K\cap \{1,2,\cdots,i-1\})$ if $i,i+1\in K$, and $j = \#(K\cap \{1,2,\cdots,i-1\})$ if $i,i+1\notin K$.
\end{enumerate}

\begin{proof}
(a) follows from the definition of $w_K$. For (b) and (c), $s_ks_{k+1}\cdots s_{n-1}t$ commutes with $s_i$ for $k>i+1$, and the claims follow from the assumptions. Finally, we have $i\neq n$ and either $i,i+1\in K$, or $i,i+1\notin K$. We shall compute $w_K^{-1}s_iw_K$. The factors starting from $s_k$ where $k>i+1$ commute with $s_i$ and do not contribute to the conjugation. For the first case, it is straightforward to check via cycle types that 
$$(s_is_{i+1}\cdots s_{n-1}t)^{-1}(s_{i+1}\cdots s_{n-1}t)^{-1}s_i(s_{i+1}\cdots s_{n-1}t)(s_is_{i+1}\cdots s_{n-1}t) = s_{n-1}.$$
Then the claim follows from the fact that 
$$(s_k\cdots s_{n-1}t)^{-1}s_j(s_k\cdots s_{n-1}t) = s_{j-1}$$
if $k<j$ and that $\#(K\cap\{1,2,\cdots,i-1\})<n-1$. The second case also follows from the last conjugation equality. 
\end{proof}
\end{lemma}

\begin{corollary}
    Each element $w_K$ is the minimal representative in its double coset $w_KW_J$, and $W/W_J = \{w_KW_J|K\subseteq\{1,2,\cdots,n\}\}$.

    \begin{proof}
        Let $Y$ be the set of minimal representatives of $W/W_J$, and $Y_k$ be the subset of $Y$ with length of element equal to $k$. We shall show by induction that $Y_k\subseteq\{w_K|K\subseteq\{1,2,\cdots,n\}\}$. Now $Y_0 = \{e\}$ and $e = w_\emptyset$. Suppose that $Y_k\subseteq\{w_K|K\subseteq\{1,2,\cdots,n\}\}$. It follows from \cite[\S2.1 Algorithm B]{GP00} that
        $$Y_{k+1}=\{sw|w\in Y_k, s\in I, l(sw) = l(w)+1\text{ and }l(sws^\prime)>l(ws^\prime),\ \forall s^\prime\in J\}.$$
        By the lemma above, $Y_{k+1}\subseteq \{w_K|K\subseteq\{1,2,\cdots,n\}\}$. Therefore, $Y\subseteq\{w_K|K\subseteq\{1,2,\cdots,n\}\}$. The equality follows from the fact that the cardinality of both sides is $2^n$. 
    \end{proof}
\end{corollary}

We obtain the set of minimal representatives in $W_J\backslash W/W_J$ as a byproduct. Let $x_0 = e$, $x_1 = t$, and $x_{i+1} = ts_{n-1}s_{n-2}\cdots s_{n-i}x_i$ for $i=1,2,\cdots,n-1$. 

\begin{corollary}
Each element $x_i$ is the minimal representative in its double coset $W_Jx_iW_J$, and $W_J\backslash W/W_J = \{W_Jx_iW_J|0\leq i\leq n\}$.

\begin{proof}
First observe that $x_k = w_{\{n-k+1,n-k+2,\cdots,n\}}$. Since $w_K$ is minimal representative in $w_KW_J$, $W_Jx_kW_J = W_Jx_lW_J$ if and only if $x_lx_k^{-1}\in W_J$, but this cannot be true unless $w_{\{n-k+1,n-k+2,\cdots,n\}} = w_{\{n-l+1,n-l+2,\cdots,n\}}$, or equivalently, $k=l$.
\par Write $K$ as a decreasing sequence $(\lambda_1>\lambda_2>\cdots>\lambda_{\#K})$. We can rewrite $w_K$ by braid relation as $\left(\prod_{j=1}^{\#K}s_{\lambda_j}\cdots s_{n-j}\right)x_{\#K}$, and $w_K\in W_Jx_{\#K}W_J$. Therefore, $W = \bigsqcup_{0\leq k\leq n}W_Jx_kW_J$.
\par Now it remains to show that $x_k$ is minimal in its double coset. It is equivalent to show that $l(s_ix_ks_j)\geq l(x_l)$ for $s_i,s_j\in J$ and $k = 0,1,..,n$. By \Cref{left_multi}, if $i=n-k$, then $s_ix_k = x_{k+1}$ and $l(s_ix_ks_j) = l(x_k)+2$ since $x_k$ is a minimal in $x_kW_J$. Otherwise, $s_ix_k = x_ks_p$ for some $s_p\in J$. Hence, either $s_ix_ls_j = x_k$ if $j = p$ and $l(s_ix_ls_j) = l(x_l)+2$ if $j\neq p$.
\end{proof}

\end{corollary}

In particular, $w_K\in W_Jx_{\#K}W_J$.

\subsection{Length of Maximal Representatives}
The length of the maximal representative in a coset is crucial in the dimension formula of some varieties considered in \Cref{decomposition}. 

The maximal representative in a single parabolic coset is exactly the product of the minimal representative and the longest element in $W_J$, yelling the following proposition.

\begin{proposition}
Let $\overline{w_K}$ be the maximal representative in the single parabolic coset $w_KW_J$. Then $l(\overline{w_K}) = (n+1)\#K-\sum_{\lambda\in K}\lambda+\frac{n(n-1)}{2}$.
\end{proposition}

Before we give a formula for the length of maximal representatives in double parabolic cosets, we introduce some notations. Let $J_1,J_2\subseteq I$ and denote $X_{J_1}$ the set of minimal representatives of $W/W_{J_1}$, $X_{J_1J_2}$ the set of minimal representatives of $W_{J_1}\backslash W/W_{J_2}$. If $J_3\subseteq J_1$, we denote $X_{J_3}^{J_1}$ the set of minimal representatives of $W_{J_1}/W_{J_3}$. We have the following Mackey decomposition for a general Coxeter group \cite[Lemma. 2.1.9.]{GP00}.

\begin{lemma}
For $d\in W$, set $J_1^d = d^{-1}J_1d$. Then $$X_{J_1} = \bigsqcup_{d\in X_{J_1J_2}}d X^{J_2}_{J_1^d\cap J_2},$$
where the multiplication is length-additive. In particular, each $w\in W$ has a unique decomposition $w = udv$ where $u\in W_{J_1}$, $d\in X_{J_1J_2}$ and $v\in X^{J_1}_{J_1^d\cap J_2}$ and $l(w) = l(u)+l(d)+l(v)$. 
\end{lemma}

Using this lemma, we can find the maximal representatives and their lengths in $W_J\backslash W/W_J$.

\begin{proposition} \label{length}
Let $\overline{x_k}$ be the longest element in $W_Jx_kW_J$. Then $$l(\overline{x_k}) = \frac{n(n-1)}{2}+\frac{k(k+1)}{2}+nk-k^2.$$

\begin{proof}
Similar to \Cref{left_multi}, one can calculate via cycle types and claim that $J\cap x_k^{-1}Jx_k = J\setminus\{s_{n-k}\}$ if $k<n$ and $J\cap x_n^{-1}Jx_n = J$. By \cite[Lemma. 2.2.1.(a)]{GP00}, the longest element in $X^J_{J\setminus\{s_{n-k}\}}$ is $(w_{J\setminus\{s_{n-k}\}})_0(w_J)_0$, and of the length $$\dfrac{n(n-1)}{2}-\dfrac{k(k-1)}{2}-\dfrac{(n-k)(n-k-1)}{2} = nk-k^2.$$
The conclusion follows from the last lemma.
\end{proof}

\end{proposition}

\subsection{An Invariant of the Minimal Representatives} \label{combinatorial_invariant}
For $w_K$ defined above, we give a formula for the function $f_K$ on $\{1,2,\cdots,n\}$ defined by $f_K(j) = \#\{i|1\leq i\leq n, 0\leq w_K(i)\leq j\}$, which is crucial for us to determine the Bruhat order and test which Schubert cell and opposite Schubert cell a Lagrangian subspace lies in. 

\begin{lemma} \label{send-to-negative}
For $i>0$, we have $w_K(i)<0$ if and only if $i+\#K >n$.

\begin{proof}
If we write $w_K$ as a product of cycles $s_ks_{k+1}\cdots s_{n-1}t =(\overline{k,k+1,\cdots,n})$ where $k\in K$, the smallest number in the rightmost cycle (the largest number in $K$) has to be strictly greater than $n-\#K$. If $i>n-\#K$, then it goes through $n-i+1$ cycles non-trivially and becomes negative, and then remains unchanged. Since there are at most $\#K$ moves, those $i$ with $i\leq n-\#K$ will never become negative.  
\end{proof}

\end{lemma}

The above analysis also gives another description of $K$ as $K = \{k|\exists i, 1\leq i\leq n\text{ and }w_K(i) = -k\}$. As a corollary, we have

\begin{corollary}
We have $w_K(\{1,2,\cdots,n\})\cap\{1,2,\cdots,n\} = \{1,2,\cdots,n\}\setminus K$.
\end{corollary}

Now we give a formula for $f_K$.

\begin{proposition}
   Let $K_{\leq j} = K\cap\{1,2,\cdots,j\}$. Then $f_K(j) = j-\#K_{\leq j}$.
   
   \begin{proof}
   We prove the statement by induction. It follows from the definition that $f_K(1) = 1-\#K_{\leq 1}$. We obverse that
   $$f_K(j+1)-f_K(j) = \left\{
   \begin{array}{c c}
      1,  & \text{if }\exists\ i>0, w_K(i) = j+1, \\
      0,  & \text{otherwise,}
   \end{array}\right.$$
and 
   $$\#K_{\leq j+1}-\#K_{\leq j} = \left\{
   \begin{array}{c c}
      1,  & \text{if }j+1\in K, \\
      0,  & \text{otherwise.}
   \end{array}\right.$$
   By the corollary above, we have that $f_K(j+1)-f_K(j)=1-(\#K_{\leq j+1}-\#K_{\leq j})$. Now the statement follows by induction.
   \end{proof}
\end{proposition}

It is easy to see from the formula that $f_K\neq f_L$ if $K\neq L$. This function can be regarded as an invariant for the cosets in $W/W_J$.

 For $K\subseteq\{1,2,..,n\}$, denote $K^\vee := \{1,2,\cdots,n\}\setminus K$ the complement of $K$. Let $w_0$ be the longest element in $W$.

\begin{lemma} \label{dual}
We have $w_0w_KW_J= w_{K^\vee}W_J$ and $W_Jw_0x_kW_J = W_Jx_{n-k}W_J$.

\begin{proof}
We need to show that $w_{K^\vee}^{-1}w_0w_K\in W_J$. It is equivalent to show that $w_{K^\vee}^{-1}w_0w_K$ sends positive number to positive number. For $i>n-\#K$, $w_K(i) = -k$ for some $k\in K$, and then $w_0(-k) = k$. Since we have that $w_{K^\vee}(\{1,2,\cdots,n\})\cap\{1,2,\cdots,n\} = K$, we have that $w^{-1}_{K^\vee}(K)>0$, so $w_{K^\vee}^{-1}w_0w_K(i)>0$.

For $i\leq n-\#K$, we have $w_K(i)\in K^\vee$. We need to show that $-w_{K^\vee}^{-1}w_K(i)>0$, or equivalently, $w_{K^\vee}^{-1}w_K(i)<0$, but this follows from the fact that $w_K(\{1,2,\cdots,n\})\cap\{1,2,\cdots,n\} = K^\vee$ and the similar equality for $K^\vee$.

The second statement follows from the first one since $x_k = w_{\{n-k+1,\cdots,n\}}$ and $w_K\in W_Jx_{\#K}W_J$. 
\end{proof}

\end{lemma}

\begin{remark}
From the formula for $f_K$, it is easy to see that $f_{K^\vee}(j) = \# K_{\leq j}$. 
\end{remark}

As an application, we can use this invariant to determine the Bruhat order of $w_K$'s. The general treatment can be found in \cite[Chapter 8]{BB05}.

\begin{corollary} \label{1-order}
We have $w_K\leq w_L$ if and only if $\#K_{\leq j}\leq\#L_{\leq j}$ for all $j=1,2,\cdots,n$, or equivalently, when both sets are written as increasing sequences, $K = (\lambda_1<\lambda_2<\cdots<\lambda_{\#K})$ and $L = (\mu_1<\mu_2<\cdots<\mu_{\#L})$, we have $\# L\geq \#K$ and $\mu_i\leq\lambda_i$ for $i = 1,2,\cdots,\#K$.
\end{corollary}

\section{Orbital decomposition of $\mathcal P_{\geq 0}$ and cellular decomposition of $U^{\pm}_{J,\geq0}$} \label{decomposition}

\subsection{Two Variations of Open Richardson Varieties}
By Bruhat decomposition, we have $$G = \bigsqcup_{w\in W/W_J}B^+\dot wP_J^+ = \bigsqcup_{w\in W_J\backslash W/W_J}P_J^+\dot wP_J^+.$$ We consider two variations of open Richardson varieties: $\tilde R_{k,l} = P_J^+\dot x_l\cdot P_J^+\cap P_J^-\dot x_k\cdot P_J^+\subseteq\mathcal P$, and $R_{K,L} = B^+\dot w_L\cdot P_J^+\cap B^-\dot w_K\cdot P_J^+$. Both of them give rise to a stratification of $\mathcal P$. It follows from the general result of open Richardson varieties that $\tilde R_{k,l}\neq\emptyset$ if and only if $x_k\leq x_l$; $R_{K,L}\neq\emptyset$ if and only if $w_K\leq w_L$. We can also describe the closure relation of these varieties in terms of Bruhat order, namely, $\tilde R_{k^\prime,l^\prime}\subseteq\overline{\tilde R_{k,l}}$ if and only if $x_k\leq x_{k^\prime}\leq x_{l^\prime}\leq x_l$, or equivalently, $k\leq k^\prime\leq l^\prime\leq l$; $R_{K^\prime,L^\prime}\subseteq\overline{R_{K,L}}$ if and only if $w_K\leq w_{K^\prime}\leq w_{L^\prime}\leq w_L$. This condition can be transformed into the combinatorial descriptions by \Cref{1-order}. 

Since $\text{Lie }P_J^++\text{Lie }P_J^- = \text{Lie }B^++\text{Lie } B^- = \text{Lie }G$, we know by Kleinman's Theorem that all the intersections are transversal, from which we can determine the dimension of these varieties.

\begin{proposition}
We have the dimension formulas
$$\dim\tilde{R}_{k,l} = nl-\frac{l(l-1)}2-\frac{k(k+1)}2,\ \dim R_{K,L} = (n+1)(\#L-\#K)+\sum_{\lambda\in K}\lambda-\sum_{\mu\in L}\mu.$$

\begin{proof}
Let $\bar{x}_k$ be the maximal element in $W_Jx_kW_J$. By Bruhat decomposition, $\dim P_J^+\dot x_l\cdot P_J^+ = l(\bar{x}_l)+\dim B^+-\dim P_J^+$. Since $P_J^- = \dot w_0^{-1}P_J^+\dot w_0$, we have $\dim P_J^-x_k\cdot P_J^+ = \dim w_0^{-1}(P_J^+\dot w_0\dot x_l\cdot P_J^+) = \dim P_J^+\dot x_{n-l}\cdot P_J^+$ by \Cref{dual}. By transversality, we have
$$\dim\tilde{R}_{k,l} = \dim P_J^+\dot x_l\cdot P_J^++\dim P_J^-\dot x_k\cdot P_J^+-\dim G,$$
and the conclusion follows from the length formula for $\bar{x}_k$, which is in \Cref{length}.

The dimension formula for $R_{K,L}$ follows from a similar argument.

\end{proof}

\end{proposition}

The purpose of this subsection is to introduce a linear algebraic criterion to test which Schubert cell a given element lies in. Let $F^{\text{std}} = (F_1^{\text{std}}\subseteq F_2^{\text{std}}\subseteq\cdots\subseteq F^{\text{std}}_n)$ be the flag of isotropic subspace where $F_i^{\text{std}}$ is spanned by $\{e_1,e_2,\cdots,e_i\}$, for $i = 1,2,\cdots,n$, and let $F^{\text{opp}}:=w_0F^{\text{std}}$ be the opposite flag, namely, $F^{\text{opp}}_i$ is spanned by $\{e_{n+1},e_{n+2},\cdots,e_{n+i}\}$ for $i = 1,2,\cdots,n$. Let $F\in\mathcal P$ be considered as a Lagrangian subspace of $\mathbb R^{2n}$.

\begin{lemma} 
 We have $F\in P^+_J\dot x_k\cdot P^+_J$ if and only if $\dim (F\cap F^{\text{std}}_n) = n-k$.

 \begin{proof}
     By \Cref{send-to-negative}, $\dim(\dot x_kF^{\text{std}}_n\cap F_n^{\text{std}}) = n-k$. If $F\in P^+_J\dot x_k\cdot F^{\text{std}}_n$, then $$\dim(F\cap F^{\text{std}}_n) = \dim(P_J^+\dot x_kF^{std}_n\cap F^{\text{std}}_n) = \dim(\dot x_kF^{\text{std}}_n\cap F^{\text{std}}_n) = n-k$$ since $F^{\text{std}}_n$ is stabilized by $P^+_J$. Since the numbers $n-k$ are different from each other, the conclusion follows.
 \end{proof}
\end{lemma}

\begin{corollary}
    We have $F\in P^-_J\dot x_k\cdot P^+_J$ if and only if $\dim (F\cap F^{\text{opp}}_n) = k$.

    \begin{proof}
       First note that $P_J^-\dot x_k\cdot P_J^+ = \dot w_0^{-1}P_J^+\dot w_0\dot x_kP_J^+ = \dot w_0^{-1}P_J^+\dot x_{n-k}P_J^+$. We have $F\in P^-_J\dot x_k\cdot P^+_J$ if and only if $\dot w_0F\in P^+_J\dot x_{n-k}\cdot P^+_J$, which is equivalent to $\dim(F\cap F^{\text{opp}}_n) = \dim(\dot w_0F\cap F_n^{\text{std}}) = k$. 
       \end{proof}
\end{corollary}

By a similar argument, we have the following criteria for $B^+\dot w_K\cdot P_J^+$ and $B^-\dot w_K\cdot P_J^+$. Recall the definition of $f_K$ from \Cref{combinatorial_invariant}.

\begin{lemma} \label{B_dim_test}
We have $F\in B^+\dot w_K\cdot P^+_J$ if and only if $\dim (F\cap F^{\text{std}}_j) = f_K(j)$ for $j=1,2,\cdots,n$.
\end{lemma}

\begin{corollary}
We have $F\in B^-\dot w_K\cdot P^+_J$ if and only if $\dim (F\cap F^{\text{std}}_j) = f_{K^\vee}(j)$ for $j=1,2,\cdots,n$.
\end{corollary}

\subsection{Orbital Decomposition}
In this section, we study the intersection of $\tilde R_{k,l}$ and $\mathcal P_{\geq 0}$. We show that they are exactly the orbits of $L_J^\circ$ on $\mathcal P_{\geq 0}$ and each of them is a connected component of the respective open Richardson variety $\tilde R_{k,l}$.

First, we recall some basic notions in linear algebra. Consider the vector space $\mathbb R^n$ equipped with a non-degenerate symmetric bilinear form $(-,-):\mathbb R^n\times\mathbb R^n\to\mathbb R$ defined by $(v,w) = \sum_{i=1}^nv_iw_i$ where $v = (v_i)_{i=1}^n$ and $w = (w_i)_{i=1}^n$. For a basis $\{a_i\}_{i=1}^n$ of $\mathbb R^n$, there exists a basis $\{a^\vee_i\}_{i=1}^n$ of $\mathbb R^n$ such that $(a_i,a^\vee_i) = 1$ and $(a_i,a^\vee_j) = 0$ whenever $i\neq j$, which we refer to as the dual basis of $\{a_i\}_{i=1}^n$. It is easy to see that if $A = [a_1,a_2,\cdots,a_n]$, then $(A^t)^{-1} = [a^\vee_1,a^\vee_2,\cdots,a^\vee_n]$.

\begin{lemma} \label{extension}
 Let $k$ and $l$ be two integers such that $0\leq k\leq l\leq n$. Given two linear independent sets $\{a_i\}_{i=1}^{n-k}$ and $\{b_i\}_{i=n-l+1}^n$ such that $(a_i,b_i) = 1$ and $(a_i,b_j)=0$ for $i\neq j$. There is a pair of basis $\{a_i\}_{i=1}^n$ and $\{b_i\}_{i=1}^n$, dual to each other, and extending $\{a_i\}_{i=1}^{n-k}$ and $\{b_i\}_{i=n-l+1}^{n}$.

\begin{proof}
For $V$ a subspace of $\mathbb R^n$, denote by $V^\perp$ the orthogonal complement of $V$. We have $\dim V+\dim V^\perp = n$. Hence
$\dim(\text{span}(b_{n-l+1},\cdots,b_n)^\perp) = l$ and $\dim(\text{span}(b_{n-l+2},\cdots,b_n)^\perp) = l-1$. There is a vector $a_{n-k+1}\in \text{span}(b_{n-l+1},\cdots,b_n)^\perp$ such that $(a_{n-k+1},b_{n-k+1}) = 1$, and $\{a_i\}_{i=1}^{n-k+1}$ is still linearly independent.

Proceeding, we may extend $\{a_i\}_{i=1}^{n-k}$ to a basis of $\mathbb R^n$ and the same procedure applies to $\{b_i\}_{i=n-l+1}^n$, completing the proof.
\end{proof}

\end{lemma}

Let $I^+_k = \text{diag}(1,1,\cdots,1,0,0,\cdots,0)$ and $I^-_k = \text{diag}(0,0,..,0,1,1,\cdots,1)$ where the number of $1$ in both matrices is $k$. Let $k$ and $l$ be two integers such that $0\leq k\leq l\leq n$. Consider $I_{k,l} = \begin{bmatrix}
    I_{n-k}^+ \\ I_l^-
\end{bmatrix}\in\mathcal P_{\geq 0}$. By the criteria from the last section, we have $I_{k,l}\in\tilde{R}_{k,l}$.

\begin{theorem} \label{orbit}
We have $L_J^\circ\cdot I_{k,l} = \mathcal P_{\geq 0}\cap\tilde{R}_{k,l}$.

\begin{proof}
First note that both of $\mathcal P_{\geq 0}$ and $\tilde{R}_{k,l}$ are invariant under the left translation of $L_J^\circ$. Hence we have that $L_J^\circ\cdot I_{k,l}\subseteq\mathcal P_{\geq 0}\cap\tilde{R}_{k,l}$.

Let $F\in\mathcal P_{\geq0}\cap\tilde{R}_{k,l}$. By the criteria from last section, $F$ can be presented by
$$\begin{bmatrix} A \\ C \end{bmatrix} = \begin{bmatrix}
    a_1 & \cdots & a_{n-l} & a_{n-l+1} & \cdots & a_{n-k} & 0 & \cdots & 0 \\
    0 & \cdots & 0 & c_{n-l+1} & \cdots & c_{n-k} & c_{n-k+1}  & \cdots & c_n
\end{bmatrix},$$
where $a_i$'s and $c_j$'s are column vectors. Since $A^tC = C^tA$, we have that $(a_i,c_j) = 0$ for $i$ and $j$ not simultaneously lying in $[n-l+1,n-k]$.

We show that the matrix $[(a_i,c_j)]_{n-l+1\leq i,j\leq n-k}$ is invertible. Suppose not. Then there is a vector $\lambda = (\lambda_i)_{i=n-l+1}^{n-k}$ such that $[(a_i,c_j)]\lambda = 0$, which implies that $\sum_{i=n-l+1}^{n-k}\lambda_ic_i$ is orthogonal to the subspace spanned by $\{a_i\}_{i=1}^{n-k}$. However, this subspace is orthogonal to $c_j$ for $j>n-k$. For the sake of dimension, $[(a_i,c_j)]$ has to be invertible.

Since $F\in\mathcal P_{\geq0}$, $[(a_i,c_j)]_{n-l+1\leq i,j\leq n-k}$ is positive definite. With the right translation by an appropriate $n\times n$ invertible matrix (equivalence in $\mathcal P$), we may further assume that $$A^tC = \text{diag}(0,\cdots,0,1,\cdots,1,0,\cdots,0),$$ or equivalently, $(a_i,c_i) = 1$ and $(a_i,c_j)=0$ for $i\neq j$. By \Cref{extension}, there is an invertible matrix $A = [a_1,a_2,\cdots,a_n]$ such that $(A^t)^{-1} = [c_1,c_2,\cdots,c_n]$. We may rearrange the order of $\{a_i\}_{i=1}^n$ and $\{c_i\}_{i=1}^n$ to assume that $\det(A)>0$. Therefore $$\begin{bmatrix} A & 0 \\ 0 & (A^t)^{-1}\end{bmatrix}\in L_J^\circ\text{, and }\begin{bmatrix} A & 0 \\ 0 & (A^t)^{-1}\end{bmatrix}\cdot I_{k,l} = F.$$ 
We conclude that $\mathcal P_{\geq 0}\cap\tilde{R}_{k,l}=L_J^\circ\cdot I_{k,l}$.

\end{proof}

\end{theorem}

We can now identify the Hausdorff closure of $\mathcal P_{>0}$ with $\mathcal P_{\geq 0}$.

\begin{corollary} \label{closure}
We have
$\overline{\mathcal P_{>0}} = \mathcal P_{\geq 0}$.

\begin{proof}
Let $V = \{(A,C)\in\mathbb R^{2n^2}|A^tC = C^tA\}$,
and $\text{Sym}(n,\mathbb R)\simeq\mathbb R^{\frac{n(n+1)}{2}}$ the space of $n\times n$ symmetric matrices. There is a right $GL_n(\mathbb R)$-action on $V$ by $(A,C).g = (Ag,Cg)$ and a right $GL_n(\mathbb R)$-action on $\text{Sym}(n,\mathbb R)$ by $S.g = g^tSg$. The map $V\to\text{Sym}(n,\mathbb R)$ sending $(A,C)$ to $A^tC$ is continuous and $GL_n(\mathbb R)$-equivariant. It induces a continuous map $\phi:\mathcal P\simeq V/\sim_{GL_n(\mathbb R)}\to\text{Sym}(n,\mathbb R)/\sim_{GL_n(\mathbb R)}$.

The projection map $\pi:\text{Sym}(n,\mathbb R)\to\text{Sym}(n,\mathbb R)/\sim_{GL_n(\mathbb R)}$, as a quotient map for a group action, is surjective and open. Since the set of positive semi-definite matrices $\text{Sym}_{\geq 0}(n,\mathbb R)$ is closed in $\text{Sym}(n,\mathbb R)$, we have that $$\pi(\text{Sym}(n,\mathbb R)\setminus\text{Sym}_{\geq 0}(n,\mathbb R))$$ is open. As its complement, $\pi(\text{Sym}_{\geq 0}(n,\mathbb R))$ is closed, and then 
$\phi^{-1}(\pi(\text{Sym}_{\geq 0}(n,\mathbb R)))=\mathcal P_{\geq 0}$ is closed, so it contains $\overline{\mathcal P_{>0}}$.

Now we show the reserve inclusion. First note that $\mathcal P_{>0}$ is $L_J^\circ$-invariant, so is $\overline{\mathcal P_{>0}}$. The conclusion follows from \Cref{orbit} if we can show that $I_{k,l}\in\overline{P_{>0}}$ for all $0\leq k\leq l\leq n$, but this follows from a direct construction of a sequence $\{I_{k,l,p}\}_{p=1}^\infty\subseteq P_{>0}$ by
$I_{k,l,p} = \begin{bmatrix}
    I_{n-k}^++\frac1p I_k^- \\ I_l^-+\frac 1p I_{n-l}^+
\end{bmatrix}.$

\end{proof}

\end{corollary}

We also have the closure relation among these orbits.

\begin{corollary}
We have 
$$\overline{\mathcal P_{\geq 0}\cap\tilde R_{k, l}} = \bigsqcup_{k\leq k^\prime\leq l^\prime\leq l}(\mathcal P_{\geq 0}\cap\tilde R_{k^\prime, l^\prime}).$$

\begin{proof}
By Bruhat decomposition, we have
$$\overline{\mathcal P_{\geq 0}\cap\tilde R_{k, l}} \subseteq  \mathcal P_{\geq 0}\cap\overline{\tilde R_{k, l}} = \bigsqcup_{k\leq k^\prime\leq l^\prime\leq l}(\mathcal P_{\geq 0}\cap\tilde R_{k^\prime, l^\prime}).$$

We note that
$$\begin{bmatrix}
I_{n-k^\prime } & & \\
 & \frac{1}{p}I_{k-k^\prime} & \\
 & & 0 \\
 0 & & \\
 & \frac{1}{p}I_{l-l^\prime} & \\
 & & I_{l^\prime}
\end{bmatrix} \in\tilde R_{k,l}, \text{ for }p = 1,2,3,\cdots.$$
The limit of this sequence is $I_{k^\prime, l^\prime}\in \mathcal P_{\geq0}\cap R_{k^\prime,l^\prime}$. Now the reverse inclusion follows from that $L_J^\circ\cdot I_{k,l} = \mathcal P_{\geq 0}\cap\tilde R_{k,l}$ for every $k\leq l$.

\end{proof}

\end{corollary}

Now we turn to show that each orbit is a connected component of the respective $\tilde R_{k,l}$. 

\begin{lemma}
The dimension of $\mathcal P_{\geq 0}\cap\tilde{R}_{k,l}$ is the same as the dimension of $\tilde{R}_{k,l}$.

\begin{proof}
We begin by determining the stabilizer of each $I_{k,l}$ in $L_J^\circ$. If $\text{diag}(A,(A^t)^{-1})$ stabilizes $I_{k,l}$, then $A$ stabilizes the subspace spanned by $\{e_1,e_2,\cdots,e_{n-k}\}$. Matrices with this property form a parabolic subgroup of $L_J^\circ$. They are of the block form 
$\begin{bmatrix}
A_1^\prime & A_2^\prime \\ 0 & A_3^\prime   
\end{bmatrix}$, where $A_1$ is a square matrix of size $n-k$. In the meantime, $(A^t)^{-1}$ stabilizes the subspace spanned by $\{e_{n-l+1},e_{n-l+2},\cdots,e_n\}$, or equivalently, $A^t$ stabilizes the same subspace. Matrices that satisfy this property are of the block form
$\begin{bmatrix}
A_1^{\prime\prime} & 0 \\ A_2^{\prime\prime} & A_3^{\prime\prime}  \end{bmatrix}$, where $A_3^{\prime\prime}$ is a square matrix of size $l$. In conclusion, $A$ has the block form
$$\begin{bmatrix}
A_1 & A_2 & A_3 \\ 0 & A_4 & A_5 \\ 0 & 0 & A_6   
\end{bmatrix},$$
where $A_1$, $A_4$, and $A_6$ are all invertible square matrices of sizes $n-l$, $l-k$, and $k$, respectively.

We show that matrices of this block form stabilize $I_{k,l}$ if and only if $A_4$ is an orthogonal matrix. By definition, $A$ stabilizes $I_{k,l}$ if and only if there is an invertible $n\times n$ matrix $B$ (written in respective block form) such that
$$\begin{bmatrix}
A_1 & A_2 & 0 \\ 0 & A_4 & 0 \\ 0 & 0 & 0   
\end{bmatrix}\begin{bmatrix}
B_1 & B_2 & B_3 \\ B_4 & B_5 & B_6 \\ B_7 & B_8 & B_9   
\end{bmatrix} = \begin{bmatrix}
I_{n-k} & 0 \\  0 & 0   
\end{bmatrix}\text{ and }
\begin{bmatrix}
0 & 0 & 0 \\ 0 & (A_4^t)^{-1} & 0\\ 0 & A_5^\prime & (A_6^t)^{-1}   
\end{bmatrix}\begin{bmatrix}
B_1 & B_2 & B_3 \\ B_4 & B_5 & B_6 \\ B_7 & B_8 & B_9   
\end{bmatrix} = \begin{bmatrix}
0 & 0 \\  0 & I_l   
\end{bmatrix}.$$
Here, $A_5^\prime$ comes from $(A^t)^{-1}$. 

Solving this matrices equations system, we get $B_3 = B_4 = B_6 = B_7 = 0$ and $B_1 = A_1^{-1}$, $B_9 = (A_6)^t$. Most importantly, $B_5 = A_4^{-1} = A_4^t$, which is equivalent to the condition that $A_4$ is orthogonal. With all the above conditions, $B_2$ and $B_8$ can be uniquely determined.

Now we can calculate that the dimension of the stabilizer of $I_{k,l}$ is
$$(n-l)^2+\frac{(l-k)^2-(l-k)}2+k^2+(n-l)(l-k)+k(n-l)+k(l-k) = n^2-nl+\frac{l(l-1)}{2}+\frac{k(k+1)}{2},$$
yelling that $\dim L_J^\circ\cdot I_{k,l} = \dim(\tilde{R}_{k,l}\cap\mathcal P_{\geq 0}) = \dim\tilde{R}_{k,l}$.

\end{proof}

\end{lemma}

\begin{proposition} \label{connected_component}
The space $\mathcal P_{\geq0}\cap\tilde{R}_{k,l}$ is a connected component of $\tilde{R}_{k,l}$.

\begin{proof}
Since $\mathcal P_{\geq0}\cap \tilde{R}_{k,l}$ is an $L^\circ_J$-orbit, it is connected. It is closed in $\tilde{R}_{k,l}$ by definition. Since it has the same dimension as $\tilde{R}_{k,l}$, it is open in $\tilde{R}_{k,l}$ by the invariance of domain.
\end{proof}

\end{proposition}

As a corollary, we obtain an orbital decomposition of $U_{J,\geq 0}^-$. It can be regarded as indexed by the rank of the positive semi-definite symmetric matrix.

\begin{corollary}
The set of $L_J^\circ$-orbit on $U^-_{J,\geq 0}$ consists of $U^-_{J,\geq 0}\cap P_J^+x_kP_J^+$ for $0\leq k\leq n$, and the space $U^-_{J,\geq 0}\cap P_J^+x_kP_J^+$ is a connected component of $U^-_J\cap P^+_Jx_kP^+_J$. We have
$$\overline{U_{J,\geq 0}^-\cap P^+_J\dot x_lP^+_J} = \bigsqcup_{k\leq l}U^-_{J,\geq0}\cap P_J^+\dot x_kP_J^+.$$
\end{corollary}

\begin{remark}
There is a similar decomposition on $U_{J,\geq 0}^+$ by intersecting $P^-_Jx_kP^-_J$.
\end{remark}

\subsection{Cellular Decomposition of $U^-_{J,\geq0}$}
It was shown in \cite[Proposition 2.7]{L94} that the totally nonnegative unipotent monoid $U_{\geq 0}^-$ admits a cellular decomposition by intersecting Bruhat cells, namely, $U_{\geq 0}^-\cap B^+\dot wB^+\simeq\mathbb R_{>0}^{l(w)}$. In this section, we prove an analog result for $U_{J,\geq 0}^-\cap B^+\dot w_KP_J^+$ by relating it to the classical Cholesky decomposition.

For $K\subseteq\{1,2,\cdots,n\}$, define
$$M_K = \{A\in \text{Mat}_{n\times n}(\mathbb R)|A_{i,j} = 0\text{ for }i<j\text{ or }j\notin K\text{ and }A_{k,k}>0\text{ for }k\in K\}.$$ Consider the map $m:\bigsqcup_{K}M_K\to\text{Sym}(n,\mathbb R)$ by $m(A) = A A^t$. The following proposition is called the Cholesky decomposition \cite[\S3.2.2]{Gen98}.

\begin{proposition}
The map $m$ is injective with image $\text{Sym}_{\geq 0}(n,\mathbb R)$.
\end{proposition}

If we endow each $M_K$ with the natural topology, then $M_K\simeq\mathbb R_{>0}^{\#K}\times\mathbb R^{\sum_{k\in K}(n-k)}$ is a topological cell, and $m$ is a homeomorphism to $\text{Sym}_{\geq 0}(n,\mathbb R)$. This gives a cellular decomposition on $\text{Sym}_{\geq 0}(n,\mathbb R)$. However, in contrast to total positivity, these cells are not semi-algebraic. 

Denote $C_K:=m(M_K)$. The main result of this section is

\begin{proposition} \label{cellular}
The intersection $U_{J,\geq 0}^-\cap B^+\dot w_KP_J^+$ is a connected component of $U_J^-\cap B^+\dot w_KP_J^+$, and
$$U_{J,\geq 0}^-\cap B^+\dot w_KP_J^+ = \left\{\begin{bmatrix}
I & 0 \\
S & I 
\end{bmatrix}\bigg|S\in C_K\right\}\simeq\mathbb R_{>0}^{\#K}\times\mathbb R^{\sum_{k\in K}(n-k)}.$$

\begin{proof}
Assume that the second statement holds. Then we can compare the dimensions and argue that $U_{J,\geq 0}^-\cap B^+\dot w_KP^+_J$ is open in $U^-_J\cap B^+w_KP^+_J$ by the invariance of domain. It is, however, also closed by definition, and thus a connected component. 

The morphism $U_J^-\to B^-\cdot P_J^+\subset\mathcal P$ by $u\to u\cdot P_J^+$ is an open embedding, via which we can translate the statement and work on the flag variety rather than the group. 

We will show that
$$\left\{\begin{bmatrix}
I  \\
S  
\end{bmatrix}\in\mathcal P\bigg|S\in C_K\right\}\subseteq U_{J,\geq 0}^-\cdot P_J^+\cap B^+\dot w_K\cdot P_J^+\text{ for }K\subseteq\{1,2,\cdots,n\}.$$ Then the two sets are identical by Cholesky decomposition. 

Recall that $F_j^{\text{std}}$ is the subspace spanned by $\{e_1,e_2,\cdots, e_j\}$. Note that, in general, if $F^\prime\in U^-_J\cdot P_J^+$ is represented by $\begin{bmatrix} I \\ S\end{bmatrix}$ with $S = [u_1,u_2,\cdots,u_n]$ when written as a collection of column vectors, then $\dim(F^\prime\cap F^{\text{std}}_j) = j-\text{rank}(u_1,u_2,\cdots,u_j)$. 

Now suppose that $F\in U_{J,\geq 0}^-\cdot P_J^+$ is represented by $\begin{bmatrix} I\\ AA^t \end{bmatrix}$ with $A\in M_K$. Writing matrix as a collection of column vectors, we have $AA^t = [\sum_{i=1}^j A_{i,j}v_i]_{j=1}^n$, where $v_i = (A_{i,1},A_{i,2},\cdots, A_{i,n})^t$. The remark above implies that
$$\dim(F\cap F^{\text{std}}_{j+1})-\dim(F\cap F^{\text{std}}_j) = \left\{\begin{array}{lr}
    1,\text{ if } A_{j,j} = 0, \\
    0, \text{ if } A_{j,j} \neq 0.
\end{array}\right.$$
By the definition of $M_K$, we have that $A_{j,j}\neq 0$ if and only if $j\in K$. By induction, $\dim(F\cap F_j^{\text{std}}) = j-\#K_{\leq j}$, and the conclusion follows from \Cref{B_dim_test}.

\end{proof}

\end{proposition}

As a corollary, we have the closure relations among these cells.

\begin{corollary}
We have $$\overline{U_{J,\geq 0}^-\cap B^+\dot w_LP_J^+} = \bigsqcup_{w_K\leq w_L}(U_{J,\geq 0}^-\cap B^+\dot w_KP_J^+).$$

\begin{proof}
By Bruhat decomposition, we have $$\overline{U_{J,\geq 0}^-\cap B^+\dot w_LP_J^+}\subseteq U_{J,\geq0}^-\cap\overline{B^+\dot w_LP_J^+} = \bigsqcup_{w_K\leq w_L}(U_{J,\geq 0}^-\cap B^+\dot w_KP_J^+).$$

Write the sets $K$ and $L$ as increasing sequences: $K = (\lambda_1<\lambda_2<\cdots<\lambda_{\#K})$ and $L = (\mu_1<\mu_2<\cdots<\mu_{\#L})$. We recall from \Cref{1-order} that $w_K\leq w_L$ is equivalent to $\# L\geq \#K$ and $\mu_i\leq\lambda_i$ for $i = 1,2,\cdots,\#K$. This defines a partial order on $2^{\{1,2,\cdots,n\}}$ generated by $K\lessdot L$ if $K = L\setminus\{i\}$ or $K = L\setminus\{i\}\cup\{i+1\}$ for $i\in L$. To show that $U_{J,\geq 0}^-\cap B^+\dot w_KP_J^+\subseteq\overline{U_{J,\geq 0}^-\cap B^+\dot w_LP_J^+}$ in general, we may assume that $K\lessdot L$.

If $K = L\setminus\{i\}$ for some $i\in L$, then the conclusion just follows from the parametrization of these cells given in \Cref{cellular}. Now suppose that $K = L\setminus\{i\}\cup\{i+1\}$ and $i+1\notin L$ (otherwise $K = L\setminus\{i\}$). For $S\in C_K$, we can find $A\in M_K$ such that $S = AA^t$. Let $D$ be the permutation matrix of $(i,i+1)$, and let $A^\prime = AD$. Then $A^\prime A^{\prime,t}$ clearly lies in $\overline{C_L}\subseteq\text{Sym}_{\geq0}(n,\mathbb R)$. But $A^\prime A^{\prime,t} = AD^2 A^t = AA^t = S$, so $C_K\subseteq\overline{C_L}$, and thus $U_{J,\geq 0}^-\cap B^+\dot w_KP_J^+\subseteq\overline{U_{J,\geq 0}^-\cap B^+\dot w_LP_J^+}$. 
 
\end{proof}

\end{corollary}

Similar to \Cref{connected_component}, we have

\begin{corollary}
The cell $U_{J,\geq0}^-\cap B^+w_KP_J^+$ is a connected component of $U_J^-\cap B^+w_KP_J^+$.
\end{corollary}

\begin{remark}
There is a similar decomposition of $U^+_{J,\geq 0}$ by intersecting $B^-\dot w_KP_J^-$.
\end{remark}

Compared to the parallel theory in total positivity, one may want to consider the decomposition $\mathcal P_{\geq 0} = \bigsqcup(\mathcal P_{\geq 0}\cap R_{K,L})$. However, a stratum on the right-hand side may not be a topological cell. For example, when $n = 2$, $K = \{2\}$ and $L = \{1\}$, we have $\mathcal P_{\geq 0}\cap R_{K,L}\simeq\mathbb R^*$. Replacing $R_{K,L}$ by the projected Richardson varieties does not help. For example, when $n=2$, the intersection of $\mathcal P_{\geq 0}$ with the projected Richardson variety of maximal dimension is isomorphic to $\mathbb R^*\times\mathbb R_{>0}$. 
It is likely that these intersections can be further decomposed as union of cells.

\section{Topology} \label{topology}
It was proven in \cite{GKL22a,GKL19} that the totally nonnegative part of a partial flag variety is homeomorphic to a closed ball. We use the techniques there to construct a contractive flow on $\mathcal P_{\geq 0}$ and show that it is homeomorphic to a closed ball. Now we review the main tool in \cite{GKL22a}.

\begin{definition}
We say a map $f:\mathbb R\times\mathbb R^N\to\mathbb R^N$ is a contractive flow if the following conditions are satisfied:

\begin{enumerate}[label = (\alph*)]
    \item the map $f$ is continuous;
    \item for all $p\in\mathbb R^N$ and $x_1,x_2\in\mathbb R$, we have $f(0,p) = p$ and $f(x_1+x_2,p) = f(x_1,f(x_2,p))$;
    \item for all $p\neq 0$ and $x>0$, we have $\Vert f(x,p)\Vert <\Vert p\Vert$.
\end{enumerate}

\end{definition}

\begin{proposition} \label{flow}
Let $Q\subseteq\mathbb R^N$ be a smooth embedded submanifold of dimension $d\leq N$, and $f:\mathbb R\times\mathbb R^N$ a contractive flow. Suppose that $Q$ is bounded and satisfies the condition
$$f(x,\bar Q)\subseteq Q,\ \forall x>0.$$
Then the closure $\bar Q$ is homeomorphic to a closed ball of dimension $d$, and $\bar Q\setminus Q$ is homeomorphic to a sphere of dimension $d-1$. 
\end{proposition}

We start to embed $\mathcal P_{\geq 0}$ to some $\mathbb R^N$ and construct the contractive flow. 
 
Let $\tau = \begin{bmatrix} 0&I\\ I&0 \end{bmatrix}\in\mathfrak{sp}_{2n}(\mathbb R)$. Being symmetric, $\tau$ has a orthogonal basis of eigenvectors
$$\{v_i:=e_{i}+e_{n+i}|1\leq i\leq n\}\cup\{v_{n+i}:=e_i-e_{n+i}|1\leq i\leq n\}$$
in $\mathbb R^{2n}$, where $v_i$ has eigenvalue $1$ and $v_{n+i}$ has eigenvalue $-1$ for $i=1,2,\cdots,n$. 

If we regard $\tau$ as a matrix in $GL_{2n}(\mathbb R)$, the left translation by $\tau$ on $Gr(n,2n)$ preserves $\mathcal P$, $\mathcal P_{\geq 0}$, and $\mathcal P_{>0}$. We shall consider the left translation of $\exp(x\tau)\in Sp_{2n}(\mathbb R)$ on $\mathcal P$.

\begin{lemma} \label{contract}
For $X\in\mathcal P_{\geq 0}$ and $x>0$, we have $\exp(x\tau)X\in\mathcal P_{>0}$.

\begin{proof}
It is straightforward to compute that 
$$\exp(x\tau) = \begin{bmatrix}\dfrac{e^x+e^{-x}}{2}I&\dfrac{e^x-e^{-x}}{2}I \\ \dfrac{e^x-e^{-x}}{2}I&\dfrac{e^x+e^{-x}}{2}I
\end{bmatrix},$$
and
$$\exp(x\tau)\begin{bmatrix}A\\ C\end{bmatrix} = \begin{bmatrix}\dfrac{e^x+e^{-x}}{2}A+\dfrac{e^x-e^{-x}}{2}C\\ \dfrac{e^x-e^{-x}}{2}A+\dfrac{e^x+e^{-x}}{2}C
\end{bmatrix}.$$

Note that $A^tC = C^tA$. We can compute that
$$\left(\dfrac{e^x+e^{-x}}{2}A+\dfrac{e^x-e^{-x}}{2}C\right)^t\left(\dfrac{e^x-e^{-x}}{2}A+\dfrac{e^x+e^{-x}}{2}C\right) = \dfrac{e^{2x}+e^{-2x}}{2}A^tC+\dfrac{e^{2x}-e^{-2x}}{4}(A^tA+C^tC).$$

By assumption, $A^tC$ is positive semi-definite. For $x>0$, each of the three terms is clearly positive semi-definite. Now it is sufficient to show that $A^tA+C^tC$ is positive definite. Since $v^tA^tAv = 0$ if and only if $Av = 0$, the claim follows from the fact that $\ker A\cap\ker C = 0$ since $(A,C)^t$ is of rank $n$. 

\end{proof}

\end{lemma}

Let $\text{Mat}_{n\times n}(\mathbb R)$ be the space of $n\times n$ matrices, and define $\psi:\text{Mat}_{n\times n}(\mathbb R)\to Gr(n,2n)$ by $\psi(B):=\text{span}(v_i+\sum_{j=1}^n B_{i,j}v_{n+i},1\leq i\leq n)$. Then the image $\psi(\text{Mat}_{n\times n}(\mathbb R))$ is the big Schubert cell being
$$\{X\in Gr(n,2n):X\cap\text{span}(v_{n+1},v_{n+2},\cdots,v_{2n}) = 0\}.$$

\begin{lemma}
The image $\psi(\text{Mat}_{n\times n}(\mathbb R))$ contains $\mathcal P_{\geq 0}$.

\begin{proof}
It is sufficient to show that for any $X\in\mathcal P_{\geq 0}$, $X\cap\text{span}(v_{n+1},v_{n+2},\cdots,v_{2n}) = 0$. Suppose $v = \sum_{i=1}^na_iv_{n+i}\in X$, then we can choose a representative $(A,C)^t$ of $X$ such that the first column of $A$ is $(a_1,a_2,\cdots,a_n)^t$ and the first column of $C$ is $(-a_1,-a_2,\cdots,-a_n)^t$, and then $(A^tC)_{1,1} = -\sum_{i=1}^na_i^2$. Since $A^tC$ is symmetric and positive semi-definite, we have  $(A^tC)_{1,1}\geq 0$, yelling that $a_i = 0$ for $i=1,2,\cdots,n$.
\end{proof}

\end{lemma}

Since $\psi$ is a smooth embedding, the restricted inverse $\psi|^{-1}_{\text{Im }\psi}$ embeds $\mathcal P_{\geq 0}$ into $\text{Mat}_{n\times n}(\mathbb R)\simeq\mathbb R^{n^2}$.

\begin{theorem} \label{ball}
The space $\mathcal P_{\geq 0}$ is homeomorphic to a closed ball of dimension $\dfrac{n(n+1)}{2}$, and $\mathcal P_{\geq 0}\setminus \mathcal P_{>0}$ is homeomorphic to a sphere of dimension $\dfrac{n^2+n-2}{2}$. 

\begin{proof}
By definition, $\exp(x\tau)(v_i+\sum_{j=1}^nB_{i,j}v_{n+j}) = e^x(v_i+e^{-2x}\sum_{j=1}^nB_{i,j}v_{n+j})$. If we define $f:\mathbb R\times \text{Mat}_{n\times n}(\mathbb R)\to \text{Mat}_{n\times n}(\mathbb R)$ by $f(x,B) = e^{-2x}B$, then $\exp(x\tau)\psi(B) = \psi(f(x,B))$. Equip $\text{Mat}_{n\times n}(\mathbb R)\simeq\mathbb R^{n^2}$ with the norm $\Vert B\Vert = \sum_{i,j}B_{i,j}^2$. It is easy to verify that $f$ defines a contractive flow on $\text{Mat}_{n\times n}(\mathbb R)$.

Let $Q = \psi^{-1}(\mathcal P_{>0})\subseteq\mathbb R^{n^2}$. Then $\bar Q = \psi^{-1}(\mathcal P_{\geq0})$. Since $\mathcal P_{\geq0}$ is a closed subset of a compact space $\mathcal P$, it is compact itself, and thus $\psi^{-1}(\mathcal P_{\geq 0})$ is bounded. By \Cref{contract}, we have
$$f(x,\bar Q)\subseteq Q,\ \forall x>0,$$
and the conclusion follows from \Cref{flow}. The dimension formula follows from the fact that $\mathcal P_{>0}$ is homeomorphic to $U^{-}_{J,>0}$.

\end{proof}

\end{theorem}

\section{Comparison with other positive structures} \label{comparison}

In this section, we compare theta positivity with other positive structures in Lagrangian Grassmannian that have been studied before.

\subsection{Total Positivity}
In this section, we show that $\mathcal P_{\geq 0}$ contains several totally nonnegative parts subject to certain choices of pinnings. 

Recall that the minimal representative $w$ of the longest element in $W/W_J$ has a reduced expression $\underline{w} = \prod_{i=1}^ns_{n+1-i}\cdots s_{n-1}t$. Deleting all the $t$ from $\underline{w}$, we obtain a reduced expression $\underline{w^\prime}$ of the longest element in $W_J$ (of type $A_{n-1}$). Let $I$ be the $2n\times 2n$ identity matrix, and recall that the negative simple root group elements are $y_i(a) = I+aE_{i+1,i}-aE_{n+i,n+i+1}$ for $1\leq i\leq n-1$ and $y_n(a) = I+aE_{2n,n}$. We consider $$U^* = \{\prod_{i=1}^ny_{n+1-i}(a_{n+1-i,i})y_{n+2-i}(a_{n+2-i,i})\cdots y_{n-1}(a_{n-1,i})y_{n}(a_{n,i})|a_{p,q}\in\mathbb R^*\text{ for all }p,q\}.$$ The subset of $U^*$ consisting of such elements where all $a_{p,q}$ are positive is the totally nonnegative unipotent cell $U^-(w)$ defined in \cite[Proposition 2.7]{L94}.

\begin{lemma} \label{matrix_form}
Let 
$$u = \prod_{i=1}^ny_{n+1-i}(a_{n+1-i,i})y_{n+2-i}(a_{n+2-i,i})\cdots y_{n-1}(a_{n-1,i})y_{n}(a_{n,i}) =\begin{bmatrix} A&0 \\ C&(A^t)^{-1}\end{bmatrix}\in U^*,$$ and 
$$u_0 = \prod_{i=1}^ny_{n-i}(a_{n-i,i+1})y_{n+1-i}(a_{n+1-i,i+1})\cdots y_{n-1}(a_{n-1,i+1}) = \begin{bmatrix} A^\prime&0 \\ 0&(A^{\prime t})^{-1}\end{bmatrix},$$ obtained by deleting all $y_n$-factors from $u$. Then we have $A = A^\prime$, and $C$ satisfies $C_{p,q}=0$ for all $p+q<n+1$.

\begin{proof}
Consider the partial products 
$$u_{(k)} = \prod_{i=1}^ky_{n+1-i}(a_{n+1-i,i})y_{n+2-i}(a_{n+2-i,i})\cdots y_{n-1}(a_{n-1,i})y_{n}(a_{n,i}) = \begin{bmatrix}   A_{(k)} & 0\\ C_{(k)} & (A^t_{(k)})^{-1}\end{bmatrix},$$ and $u_{(k),0} = \prod_{i=1}^ky_{n-i}(a_{n-i,i+1})y_{n+1-i}(a_{n+1-i,i+1})\cdots y_{n-1}(a_{n-1,i+1})$ for $k = 1,2,\cdots,n$. We will prove by induction that $u_{(k)}$ and $u_{(k),0}$ have the same upper left block, and $C_{(k),p,q} = 0$ for $p+q<n+1$ or $q<n-k+1$.

When $k = 1$, $u_{(1)} = y_n(a_{n,1})$, and this is obvious. Now we consider

$$y_{n-k}(a_{n-k,k+1})y_{n+1-k}(a_{n+1-k,k+1})\cdots y_{n-1}(a_{n-1,k+1}) = \begin{bmatrix}   \tilde A & 0 \\ 0 & (\tilde A^t)^{-1} \end{bmatrix}.$$

We have that $\tilde A$ is unipotent lower triangular, $\tilde A_{p,q} = 0$ if $q<n-k$ and $p>q$, and that

$$u_{(k+1)} =\begin{bmatrix} A_{(k)} & 0 \\ C_{(k)} & (A^t_{(k)})^{-1} \end{bmatrix} \begin{bmatrix}    \tilde A & 0 \\ 0 & (\tilde A^t)^{-1} \end{bmatrix} \begin{bmatrix} I & 0 \\ a_{n,k+1}E_{n,n} & I \end{bmatrix} = \begin{bmatrix} A_{(k)}\tilde A & 0 \\ C_{(k)}\tilde A+a_{n,k+1}(A_{(k)}^t\tilde A^t)^{-1}E_{n,n} & (A_{(k)}^t\tilde A^t)^{-1} \end{bmatrix}.$$

Now the assertion follows from the induction assumption.
\end{proof}

\end{lemma}

A family of regular functions is defined for an arbitrary connected reductive group and studied in \cite{MR04}. In classical types, some of these functions are matrix minors. We shall use the result there to calculate some minors of matrices in $U^*$. These functions also appear in \Cref{generalized_plucker}.

First, we need some combinatorial information from the reduced expression. Recall that for a reduced expression $\underline{v} = s_{i_1}s_{i_2}\cdots s_{i_k}$, one can attach an order sequence of positive co-roots $(s_{i_2}\cdots s_{i_k}.\alpha_{i_1}^\vee,s_{i_3}\cdots s_{i_k}.\alpha_{i_2}^\vee,\cdots,\beta^\vee_{k-1} =s_{i_k}\alpha_{i_{k-1}}^\vee, \alpha_{i_k}^\vee)$, which we refer to as the reflection sequence of $\underline{v}$. The underlying set of this sequence is $\{\alpha^\vee\in(\Phi^\vee)^+|v.\alpha\in(\Phi^\vee)^-\}$.

The root systems $\Phi_J$ and $\Phi^\vee_J$ attached to $J$ are both of type $A_{n-1}$. While $\Phi$ is of type $C_n$, $\Phi^\vee$ is of type $B_n$ with the simple short root $\alpha_n^\vee$.  We denote the reflection sequence of $\underline w$ by $(\beta_i^\vee)$ and the reflection sequence of $\underline{w^\prime}$ by $(\gamma_i^\vee)$. The following lemma can be checked directly.

\begin{lemma} \label{refl_seq}
The reflection sequences for $\underline{w}$ and $\underline{w^\prime}$ are determined by
\begin{itemize}
    \item $\beta^\vee_{\frac{i(i-1)}{2}+j} = \alpha^\vee_j+\alpha^\vee_{j+1}+\cdots+\alpha^\vee_{i-1}+2\alpha^\vee_{i}+2\alpha^\vee_{i+1}+\cdots+2\alpha_n^\vee$ for $1\leq i\leq n$ and $1\leq j\leq i-1$, and $\beta^\vee_{\frac{i(i+1)}{2}} = \alpha^\vee_i+\alpha^\vee_{i+1}+\cdots+\alpha_n^\vee$;
    \item $\gamma_{\frac{i(i-1)}{2}+j}^\vee=\alpha_j^\vee+\alpha_{j+1}^\vee+\cdots+\alpha_i^\vee$ for $1\leq i\leq n-1$ and $1\leq j\leq i$.
\end{itemize}
\end{lemma}

For a matrix $M$, two subsets $K,L$ of $\{1,2,\cdots,n\}$ with $\#K = \#L$, we denote by $\triangle_{K,L}(M)$ the minors of $M$ with rows' labels taken from $K$ and columns' labels taken from $L$. 

\begin{lemma}
Retain the notations of \Cref{matrix_form}. We have the following formulas:
\begin{itemize}
    \item $$\triangle_{\{k,k+1,\cdots,n\},\{1,2,\cdots,n-k+1\}}(C) = \prod_{1\leq p\leq k}a_{n,p}\prod_{1< i\leq n,j\leq k< i}a_{n-i+j,i}\prod_{1\leq j<i\leq k}a_{n-i+j,i}^2$$ for $k = 1,2,\cdots,n$;
    \item $$\triangle_{\{k,k+1,\cdots,n\},\{1,2,\cdots,n-k+1\}}(A) = \prod_{1<i\leq n,j\leq k<i}a_{n-i+j,i}$$ for $k=1,2,\cdots,n-1$;
    \item $$\triangle_{\{1,2,\cdots,k\},\{1,2,\cdots,k\}}(A^tC) = \prod_{1\leq p\leq k}a_{n,p}\prod_{1\leq q<n, 1<i\leq n, q+i+k\geq n}a_{q,i}^2$$ for $k=1,2,\cdots,n$.
\end{itemize}

\begin{proof}
Denote by $\omega_k$ the $k$-th fundamental weights of $G$. It follows from \cite[Lemma 7.5 (1)]{MR04} that

$$\triangle_{\{k,k+1,\cdots,n\},\{1,2,\cdots,n-k+1\}}(C) = \prod_{1\leq i\leq n, 1\leq j\leq i}a_{n-i+j,i}^{\langle\omega_k,\beta^\vee_{\frac{i(i-1)}{2}+j}\rangle} = \prod_{1\leq p\leq k}a_{n,p}\prod_{1< i\leq n,j\leq k< i}a_{n-i+j,i}\prod_{1\leq j<i\leq k}a_{n-i+j,i}^2.$$

$$\triangle_{\{k,k+1,\cdots,n\},\{1,2,\cdots,n-k+1\}}(A) = \prod_{1\leq i\leq n-1, 1\leq j\leq i}a_{n-i+j-1,i+1}^{\langle\omega_k,\gamma^\vee_{\frac{i(i-1)}{2}+j}\rangle} = \prod_{1<i\leq n,j\leq k<i}a_{n-i+j,i}.$$

Here, we use \Cref{refl_seq}. When $k = n$, the minor of $A$ considered here is the determinant $1$.

Since $A$ is upper triangular and $C_{p,q}=0$ for $p+q<n+1$, we have
$$\triangle_{\{1,2,\cdots,k\},\{1,2,\cdots,k\}}(A^tC) = \triangle_{\{k,k+1,\cdots,n\},\{1,2,\cdots,n-k+1\}}(A)\triangle_{\{k,k+1,\cdots,n\},\{1,2,\cdots,n-k+1\}}(C),$$
and then

$$\triangle_{\{1,2,\cdots,k\},\{1,2,\cdots,k\}}(A^tC) = \prod_{1\leq p\leq k}a_{n,p}\prod_{1\leq q<n, 1<i\leq n, q+i+k\geq n}a_{q,i}^2.$$

\end{proof}

\end{lemma}

If we take the subset $U^*_{(n,>0)}$ of $U^*$ defined by
$$U^* = \{\prod_{i=1}^ny_{n+1-i}(a_{n+1-i,i})y_{n+2-i}(a_{n+2-i,i})\cdots y_{n-1}(a_{n-1,i})y_{n}(a_{n,i})|a_{n,q}\in\mathbb R_{>0}\text{ for all }q\}.$$

Then for every $u\in U^*_{(n,>0)}$, the respective $A^tC$ is positive definite, and it implies the following theorem.

\begin{theorem} \label{dense}
We have that $$U_{(n,>0)}^*\cdot P_J^+ \subseteq\mathcal P_{>0}.$$
\end{theorem}

The reverse inclusion does not hold. For $n=2$, a symmetric positive definite matrix obtained as $A^tC$ for such an element $u$ has an extra condition that no entry vanishes. 

Now we can compare $\mathcal P_{>0}$ with total positivity. For $\epsilon\in\{\pm 1\}$, define $y_{i,\epsilon}:\mathbb R\to G$ by $y_{i,\epsilon}(a) = y_i(\epsilon a)$ and similarly, $x_{i,\epsilon}$. Consider a collection of pinnings 
$$\mathrm C = \{(T,B^+,B^-,x_{i,\epsilon_i},y_{i,\epsilon_i};1\leq i\leq n)|\epsilon_i\in\{\pm1\}\text{ for }i\neq n, \epsilon_n = 1\}.$$
For $\mathbf P = (T,B^+,B^-,x_{i,\epsilon_i},y_{i,\epsilon_i};1\leq i\leq n)\in\mathrm C$, let $U^-(w)_{(\mathbf P,>0)}$ be the totally nonnegative unipotent cell of $w$ defined by $\mathbf P$, and let $\mathcal P_{(\mathbf P,>0)} := U^-(w)_{(\mathbf P,>0)}\cdot P_J^+$ be the totally positive part of $\mathcal P$ defined by $\mathbf P$. 

\begin{corollary} \label{TP}
We have that 
$$\bigsqcup_{\mathbf P\in\mathrm C}\mathcal P_{(\mathbf P,>0)}\subsetneqq\mathcal P_{>0},\ 
\bigsqcup_{\mathbf P\in\mathrm C}\mathcal P_{(\mathbf P,\geq0)}\subseteq\mathcal P_{\geq0}.$$
The equality of nonnegative parts holds when $n=2$, but fails when $n\geq 3$.
\end{corollary}

When $n=2$, $\bigsqcup_{\mathbf P\in\mathrm C}\mathcal P_{(\mathbf P,>0)}$ is dense in $\mathcal P_{>0}$, and thus $\bigsqcup_{\mathbf P\in\mathrm C}\mathcal P_{(\mathbf P,\geq0)}=\mathcal P_{\geq0}$. However, this is not the case in general. For example, when $n=3$, elements of the form
$$y_3(>0)y_2(>0)y_3(>0)y_1(>0)y_2(<0)y_3(>0)\cdot P_J^+$$
does not lie in the closure of the former one, as the two $y_2$-factors have different signs. For $n>2$, taking union of totally nonnegative parts subject to different pinnings will not produce $\mathcal P_{\geq 0}$.

\subsection{Generalized Pl\"ucker positivity}
Consider the natural representation $\mathbb R^{2n}$ of $Sp_{2n}(\mathbb R)$, the fundamental representation corresponding to the long simple root can be identified as the invariant subspace $V$ of $\bigwedge^n\mathbb R^{2n}$ with the highest weight vector $e_1\wedge e_2\wedge\cdots\wedge e_n$. Recall that we have a lift $\dot w$ for each $w\in W$, such that 
$$\dot{s}_i = I_{2n}-E_{i,i}-E_{i+1,i+1}-E_{i,i+1}+E_{i+1,i}-E_{n+i,n+i}-E_{n+i+1,n+i+1}-E_{n+i,n+i+1}+E_{n+i+1,n+i}, \text{ for }1\leq i\leq n-1,$$
$$\dot{s}_n = I_{2n}-E_{n,n}-E_{2n,2n}+E_{2n,n}-E_{n,2n},$$
and $\dot w = \dot w_1\dot w_2$, if $w = w_1w_2$ and $l(w) = l(w_1)+l(w_2)$. 

Let $$\underline{I} = \{(k_1,\cdots,k_n)|1\leq k_j\leq 2n,k_l<k_{l+1},k_p\neq k_q+n,\ \forall 1\leq l\leq n-1, p\neq q\},$$ and $d:\{1,2,\cdots,2n\}\to\{1,2,\cdots,n\}$ defined by 
$$d(k) = \left\{\begin{array}{lr}
    k,\text{ if } k\leq n, \\
    k-n, \text{ if } k>n.
\end{array}\right.$$
For $\underline{i}\in\underline{I}$, define $\text{sgn}(\underline{i}) = (-1)^{\#\{(p,q)|p<q,d(k_p)>d(k_q)\}}$. The set $$\{\dot w\eta_{\omega_n}|w\in W^J\} = \{\text{sgn}(\underline i)(\wedge_{k\in\underline i}e_k)|\underline i\in\underline I\}$$ is a set of extremal weight vectors. We show

\begin{proposition} \label{Plucker}
Let $\mathcal P_{\triangle>0}$ and $\mathcal P_{\triangle\geq 0}$ be the generalized Pl\"ucker positive and nonnegative parts of $\mathcal P$, respectively. Then $\mathcal P_{\triangle>0} = \mathcal P_{>0}$, and $\mathcal P_{\triangle\geq 0} = \mathcal P_{\triangle\geq 0}$.

\begin{proof}

For $F\in\mathcal P_{\triangle>0}$, $F$ can be represented as $\begin{bmatrix}
    I\\ S
\end{bmatrix}$, where $S$ is symmetric. The positivity condition now translates into the condition that $S$ is positive definite, and we have that $\mathcal P_{\triangle>0} = \mathcal P_{>0}$. The last equality is obtained by taking closure on both sides of the first equality, applying \Cref{plucker_closure} and \Cref{closure}.
    
\end{proof}

\end{proposition}

As a corollary, we can compare the total positivity and the generalized Pl\"ucker positivity on $\mathcal P$. 

\begin{corollary}
We have that 
$$\bigsqcup_{\mathbf P\in\mathrm C}\mathcal P_{(\mathbf P,>0)}\subsetneqq\mathcal P_{\triangle>0},
\ \bigsqcup_{\mathbf P\in\mathrm C}\mathcal P_{(\mathbf P,\geq0)}\subseteq\mathcal P_{\triangle\geq0}.$$
The equality of nonnegative parts holds when $n=2$, but fails when $n\geq 3$.
\end{corollary}

We point out that the Gelfand-Serganova strata do not give rise to a cellular decomposition of $\mathcal P_{\geq 0}$. When $n=2$, the intersection of $\mathcal P_{\geq 0}$ with the maximal Gelfand-Serganova stratum is homeomorphic to $\mathbb R^*\times\mathbb R_{>0}$.

\begin{remark}
Karpman \cite[Proposition 6.1]{Km18} showed that for a choice of a symplectic form and a pinning $\mathbf P$, one can identify $\mathcal P_{(\mathbf P,\geq0)}$ as the intersection of $\mathcal P$ with the totally nonnegative Grassmannian $Gr_{\geq0}(n,2n)$. Up to the choice of symplectic form, we can say that $\mathcal P_{\geq0}$ has less nonnegativity conditions than the one considered by Karpman.
\end{remark}

\subsection{The Space of Cactus Networks}
Cactus networks were introduced and studied as a generalization of electrical networks by Lam \cite{Lam18}, and he identified the space $\mathcal X_{n+1,\geq0}$ of cactus networks as a subset of the totally nonnegative Grassmannian $Gr_{\geq0}(n+2,2n+2)$. Let $\mathcal X_{n+1}$ be the Zariski closure of $\mathcal X_{n+1,\geq0}$ in $Gr(n+2,2n+2)$. It was shown by Chepuri, George, and Speyer \cite{CGS21}, and by Bychkov, Gorbounov, Kazakov, and Talalaev \cite{BGKT23}, that $\mathcal X_{n+1}\simeq Sp_{2n}/P_J^+$. Gao, Lam, and Xu \cite[\S6]{GLX22} showed that $\mathcal X_{n+1,\geq0}$ can be defined by nonnegativity of coordinates with respect to a certain basis in $V_{\omega_n}$, the fundamental representation of $Sp_{2n}$ corresponding to the long simple root. As a comparison, the theta nonnegative part $\mathcal P_{\geq0}$ is defined by nonnegativity of coordiantes with respect to extremal weight vectors, in the same representation.

\medskip

\printbibliography

\end{document}